\documentclass[11pt,a4paper,equation]{article}
\usepackage[greek,american]{babel}
\usepackage[iso-8859-7]{inputenc}
\usepackage{amssymb,amsfonts}
\usepackage[dvips]{graphicx}
\usepackage{amsmath}
\usepackage{epsfig}
\usepackage{latexsym}
\usepackage{makeidx}
\usepackage{pst-math,pst-xkey}
\newtheorem{lemma}{Lemma}[section]
\newtheorem{theorem}{Theorem}[section]

\numberwithin{equation}{section}
\makeindex
\begin{document}
\date{}
\author{ Aristides V. Doumas$^{1}$ and Vassilis G. Papanicolaou$^{2}$ \\
%EndAName
Department of Mathematics\\
National Technical University of Athens\\
Zografou Campus\\
157 80 Athens, GREECE\\
$^{1}$aris.doumas@hotmail.com \quad $^{2}$papanico@math.ntua.gr}
\title{The Coupon Collector's Problem Revisited: Generalizing the Double Dixie Cup Problem of Newman and Shepp}
\maketitle
\begin{abstract}
The ``double Dixie cup problem" of D.J. Newman and L. Shepp \cite{N-S} is a well-known variant of the coupon collector's problem, where
the object of study is the number $T_{m}(N)$ of coupons that a collector has to buy in order to complete $m$ sets of all $N$ existing different
coupons. More precisely, the problem is to determine the asymptotics of the expectation (and the variance) of $T_{m}(N)$, as well as its limit
distribution, as the number $N$ of different coupons becomes arbitrarily large. The classical case of the problem, namely the case of equal coupon
probabilities, is here extended to the general case, where the probabilities of the selected coupons are unequal.

In the beginning of the article we give a brief review of the formulas for the moments and the moment generating function of the random variable $T_{m}(N)$. Then,
we develop techniques of computing the asymptotics of the first and the second moment of $T_{m}(N)$ (our techniques apply to the higher moments of
$T_{m}(N)$ as well). From these asymptotic formulas we obtain the leading behavior of the variance $V[\,T_{m}(N)\,]$ as $N \to \infty$.
Finally, based on the asymptotics of $E[\,T_{m}(N)\,]$ and $V[\,T_{m}(N)\,]$ we obtain the limit distribution of the random variable
$T_{m}(N)$ for large classes of coupon probabilities. As it turns out, in many cases, albeit not always, $T_{m}(N)$ (appropriately normalized) converges in
distribution to a Gumbel random variable. Our results on the limit distribution of $T_{m}(N)$ generalize a well-known result of P. Erd\H{o}s
and A. R\'{e}nyi \cite{E-R} regarding the limit distribution of $T_{m}(N)$ for the case of equal coupon probabilities.
\end{abstract}

\textbf{Keywords.} Urn problems; coupon collector's problem; double Dixie cup problem; limit distribution; Gumbel
distribution; generalized Zipf law.\\\\
\textbf{2010 AMS Mathematics Classification.} 60F05; 60F99.

\section{Introduction}

The \textbf{``coupon collector's problem'' (CCP)} pertains to a population whose members are of $N$ different \emph{types} (e.g., baseball cards, viruses, fish, words, etc). For
$1 \leq j \leq N$ we denote by $p_j$ the probability that a member of the population is of type $j$, where $p_j > 0$ and $\sum_{j=1}^{N}p_{j}=1$.
We refer to the $p_j$'s as the \textit{coupon probabilities}. The members of the population are sampled independently \textit{with replacement} and their types are recorded. Naturally, a quantity of interest here is the number $T(N)$ of trials needed until all $N$ types are detected (at least once).\\
CCP belongs to the family of \textit{Urn problems\textit} among with other classical problems, such as the birthday and the matching problem.
In its simplest form (i.e. when all $p_{j}$'s are equal and the collector aims for one complete set of coupons) the problem has appeared in many standard probability textbooks
(e.g., in W. Feller's classical work \cite{F}, as well as in \cite{D}, \cite{Ma}, \cite{Ro1}, and \cite{Ro2}, to name a few).
Its origin can be traced back to De Moivre's treatise
\textit{De Mensura Sortis} of 1712 (see, e.g., \cite{Ho}) and Laplace's pioneering work \textit{Theorie Analytique de Probabilites} of 1812 (see \cite{D-H}). The problem was related to the Dixie Cup Company, since in the 1930's the
company introduced a highly successful procedure by which children collected Dixie lids to receive ``Premiums'', beginning with illustrations of their favored Dixie Circus characters, and then Hollywood stars and Major League baseball players (see \cite{Ma}, \cite{DC}).\\
CCP has attracted the attention of various researchers due to its applications to several areas of science (computer science--search algorithms, mathematical programming, optimization, learning processes, engineering, ecology, as well as linguistics---see, e.g., \cite{B-H}).

For the asymptotics of the moments, as well as for the limit distribution of the random variable $T(N)$,
there is a plethora of articles obtaining a variety of results (for the case of equal probabilities see for instance \cite{E-R}, \cite{KAP}, \cite{FLA};
as for the case of unequal probabilities, see, e.g., \cite{B} ---which presents the results of R.K. Brayton's Ph.D. thesis under the supervision of Norman Levinson--- \cite{N}, \cite{DP}, \cite{DPM},
and the references therein). \\
Several variants of the original problem have been studied. Among them there is the so-called
\textbf{``double Dixie cup problem"}, which reads: How long does it take to obtain $m$ complete sets of $N$ coupons?

Let $T_{m}(N)$ be the number of trials a collector needs in order to accomplish this goal (obviously, $T_1(N) = T(N)$, thus
the case $m=1$ reduces to the more ``classical" CCP).\\
Naturally, the simplest case occurs when one takes $p_{1}=\cdot \cdot \cdot =p_{N}=1/N$. For this case D.J. Newman and L. Shepp \cite{N-S} obtained that, for fixed $m$,
\begin{equation}
E\left[\, T_m(N)\,\right]= N \ln N + \left(m-1\right) N \ln \ln N + N C_m + o(N)
\label{1}
\end{equation}
as $N\rightarrow \infty$, where $C_{m}$ is a constant depending on $m$. Roughly speaking, formula (\ref{1}) tells us that, on the average, the first set ``costs''
$N \ln N + O(N)$, while each additional set has an additional cost of $N \ln \ln N  + O(N)$.\\
Soonafter, P. Erd\H{o}s and A. R\'{e}nyi \cite{E-R} went a step further and determined the limit distribution of $T_{m}(N)$, as well as the exact value of the constant $C_{m}$. They proved that
\begin{equation}
C_m = \gamma - \ln \left(m-1\right)!,
\label{2}
\end{equation}
where $\gamma=0.5772\cdots$ is the Euler-Mascheroni constant, and that for every real $y$ the following limiting result holds:
\begin{equation}
\lim_{N \rightarrow \infty} P\left\{\frac{T_m(N) - N \ln N - (m - 1) N \ln \ln N}{N} \leq y\right\}
= \exp \left(-\frac{e^{-y}}{(m-1)!}\right)
\label{3}
\end{equation}
or, equivalently,
\begin{equation}
\lim_{N \rightarrow \infty} P\left\{\frac{T_m(N) - N \ln N - (m - 1) N \ln \ln N + N \ln (m-1)!}{N} \leq y\right\}
= e^{-e^{-y}}
\label{333}
\end{equation}
(in the right-hand side of \eqref{333} we have the standard Gumbel distribution function; recall that its expectation is
$\gamma$ and its variance is $\pi^2 / 6$).

Later, and as long as the coupon probabilities remained equal, this result was generalized in \cite{KAP} and \cite{FLA}.

In the present paper we extend the results of Newman-Shepp \cite{N-S} and Erd\H{o}s-R\'{e}nyi \cite{E-R} to large families of unequal
coupon probabilities. Notice that in practically all applications the coupon probabilities are not equal (for example, in several applications
the coupon probabilities follow a generalized Zipf law, a case which is covered by our results). As we will see, for many families of
coupon probabilities the quantity $T_m(N)$, after an appropriate normalization, converges in distribution to the standard Gumbel random
variable as $N \to \infty$. The correct normalization of $T_m(N)$, which depends on the coupon probabilities, is determined with the help of
the asymptotics of $E[\,T_m(N)\,]$ and $V[\,T_m(N)\,]$. We also present
families of coupon probabilities for which the limit distribution of $T_m(N)$ (again after a suitable normalization) is not Gumbel.

\subsection{Moments and the moment generating function of $T_m(N)$}

Suppose that, for $j = 1, \dots, N$, we denote by $W_j$ the number of trials needed in order to obtain $m$ times the coupon of type $j$.
Then, it is clear that $W_j$ is a negative binomial random variable (with parameters $m$ and $p_j$) and
\begin{equation*}
T_m(N) = \max_{1 \leq j \leq N} W_j.
\end{equation*}
However, the above formula for $T_m(N)$ is not very useful, since the $W_j$'s are not independent. Instead, one can apply a clever
``Poissonization technique" found in \cite{Ro2} in order to get explicit formulas for the moments of $T_m(N)$.

Let $Z(t)$, $t \geq 0$, be a Poisson process with rate $\lambda = 1$. We imagine that each Poisson event associated to $Z$ is a collected coupon,
so that $Z(t)$ is the number of collected coupons at time $t$. Next, for $j = 1, \dots, N$, let $Z_j(t)$ be the number of type-$j$ coupons
collected at time $t$. Then, the processes $\{Z_j(t)\}_{t \geq 0}$, $j = 1, \dots, N$, are independent Poisson processes with rates $p_j$
respectively \cite{Ro2} and, of course, $Z(t) = Z_1(t) + \cdots + Z_N(t)$. If $X_j$, $j = 1, \dots, N$, denotes the time of the $m$-th
event of the process $Z_j$, then $X_1, \dots, X_N$ are obviously independent (being associated to independent processes) and
\begin{equation}
X := \max_{1 \leq j \leq N} X_j
\label{P1}
\end{equation}
is the time when all different coupons have arrived at least $m$ times.

Now, for each $j = 1, \dots, N$, $X_j$ is Erlang with parameters $m$ and $p_j$, meaning that
\begin{equation}
P\{X_j > t\} = S_m(p_j t) e^{-p_j t},
\label{P2}
\end{equation}
where
$S_{m}(y)$ denotes the $m$-th partial sum of $e^{y}$, namely
\begin{equation}
S_{m}(y) := 1+y+\frac{y^{2}}{2!}+\cdots+\frac{y^{m-1}}{\left(m-1\right)!}=\sum_{l=0}^{m-1}\frac{y^l}{l!}\,.
\label{7}
\end{equation}
Incidentally, let us observe that
\begin{equation}
0 < e^{-y} S_{m}(y) < 1
\qquad \text{for all }\; y > 0.
\label{7aa}
\end{equation}
It follows from (\ref{P1}) and the independence of the $X_j$'s that
\begin{equation}
P\{X \leq t\} = \prod_{j=1}^N \left[ 1 - S_m(p_j t) e^{-p_j t} \right].
\label{P3}
\end{equation}
It remains to relate $X$ and $T_m(N)$. Clearly,
\begin{equation}
X = \sum_{k=1}^{T_m(N)} U_k,
\label{P4}
\end{equation}
where $U_1, U_2, \dots$ are the interarrival times of the process $Z$. It is common knowledge that the $U_j$'s are independent and exponentially
distributed random variables with parameter $1$. In order to compute the moments of $T_m(N)$ via formula (\ref{P4}) we need the formula
\begin{equation}
E\left[ \left( \sum_{k=1}^M U_k \right)^r \, \right] = M (M+1) \cdots (M + r - 1) =: M^{(r)},
\qquad
r = 1, 2, \dots
\label{P5}
\end{equation}
(the justification of (\ref{P5}) is immediate, if we just notice that $U_1 + \cdots + U_M$ is Erlang with parameters $M$ and $1$). Since
$T_m(N)$ is independent of the $U_j$'s, formulas (\ref{P4}) and (\ref{P5}) imply
\begin{equation}
E\left[ X^r \, | \, T_m(N)\right] = T_m(N)^{(r)} E[U_j] = T_m(N)^{(r)},
\qquad
r = 1, 2, \dots,
\label{P6}
\end{equation}
hence, by (\ref{P6}), with the help of (\ref{P2}), we obtain
\begin{equation}
E\left[ T_m(N)^{(r)} \right] = E\left[ X^r \right]
= r\int_0^{\infty} \left\{1-\prod_{j=1}^{N}\left[1 - S_{m}(p_j t) e^{-p_j t} \right]\right\} t^{r-1} dt
\label{P7}
\end{equation}
for $r = 1, 2, \dots \,$. The quantity $E\left[ T_m(N)^{(r)} \right]$ is, actually, the $r$\textit{-th rising moment of} $T_m(N)$.
In particular,
\begin{align}
E[T_m(N)]&=\int_{0}^{\infty}\left\{1-\prod_{j=1}^{N}\left[1-S_{m}(p_{j}t) e^{-p_{j}t} \right]\right\} dt,
\label{5}
\\
E\left[T_m(N)\left(T_{m}(N)+1\right)\right]&=2\int_{0}^{\infty}
\left\{1-\prod_{j=1}^{N}\left[1-S_{m}(p_{j}t) e^{-p_{j}t}\right]\right\} t dt,
\label{5A}
\end{align}
and, of course,
\begin{equation}
V\left[\, T_m(N)\,\right]=E\left[\,T_m(N)\left(T_{m}(N)+1\right)\,\right]-E\left[\, T_m(N)\,\right]-E\left[\, T_m(N)\,\right]^{2}.
\label{VAR}
\end{equation}
Formulas (\ref{5}) and (\ref{5A}) were first derived in \cite{B} by a more complicated argument. As far as we know, the more general formula
(\ref{P7}) is new.

Using (\ref{P7}) one can easily obtain the moment generating function of $T_m(N)$:
\begin{equation}
G(z) := E\left[ z^{-T_m(N)} \right]
= 1 - (z-1) \int_0^{\infty} \left\{1-\prod_{j=1}^{N}\left[1 - S_{m}(p_j t) e^{-p_j t} \right]\right\} e^{-(z-1) t} dt,
\label{P8}
\end{equation}
where $\Re(z) > 1$.

\subsection{Discussion}
Under the quite restrictive assumption of ``nearly equal coupon probabilities", namely
\begin{equation}
\lambda(N) :=
\frac{\max_{1\leq j\leq N}{\left\{p_{j}\right\}}}{\min_{1\leq j\leq N}{\left\{p_{j}\right\}}}\leq M < \infty, \qquad\text{independently of $N$,}
\label{4}
\end{equation}
R.K. Brayton \cite{B} obtained detailed asymptotics of the expectation $E[T_m(N)]$, while for the asymptotics of the variance, he only did the
case $m=1$, where he found the formula
\begin{equation}
V\left[\, T_1(N)\,\right]
= N^{2}\left[\frac{\pi^{2}}{6} + O\left(\frac{\ln \ln \ln N}{\ln \ln N}\right)\right] \quad \text{as}\quad N\rightarrow \infty.
\label{4aa}
\end{equation}
The present paper builds on \cite{DP}, where the case $m=1$ was considered. Our results here are valid for
all positive integers $m$, including $m=1$.

The rest of our work is organized as follows.  In Section 2 we calculate the asymptotics of $E[T_m(N)]$ and $V[T_m(N)]$ for the  the general case of
unequal probabilities. We first show how to create a sequence $\pi_N = (p_1, \dots, p_N)$, $N = 2, 3, \dots$, of probability measures (i.e. of coupon probabilities) by successive normalizations of the
terms of a given, albeit arbitrary, sequence $\alpha = \{a_j\}_{j=1}^{\infty}$ of positive real numbers. Thus, we need to focus on the sequence
$\alpha$. First (Case I) we consider a large class of sequences $\alpha$, such that
$a_j \rightarrow \infty$. The main result for this case is presented in Theorem 2.2. For instance, polynomial and exponential families of coupon probabilities fall in this category (e.g., the so-called \textit{linear} case).

Then (Case II) we consider classes of decaying sequences $\alpha$ such that $a_j \rightarrow 0$. This case is much more challenging. It turns out that in order to obtain the leading term of the variance $V[T_{m}(N)]$ (see Theorem 2.5) one has to go deep in the asymptotics of $E\left[\, T_m(N)\,\right]$ (up to the fifth term) and $E\left[\,T_m(N)\left(T_{m}(N)+1\right)\,\right]$ (up to the sixth term). These asymptotic formulas are presented in Theorems 2.3 and 2.4 respectively. It is notable that the
\textit{generalized Zipf law} falls in this category.

The approach presented in Section 2 can be used to calculate the asymptotics, as $N \to \infty$, of the $r$-th rising moment of $T_{m}\left(N\right)$, for any positive integer $r$.

Section 3 is a short section where we present some illustrative examples. Then,
in Section 4 we take advantage of our formulas in order to find the limit distribution of $T_{m}(N)$ (appropriately normalized) for a very
large class of coupon probabilities. More precisely, for sequences of Case I the limit distribution is obtained (in Theorem 4.1) by using formula
(\ref{P8}). As for sequences of Case II, we combine our asymptotic formulas with a limit theorem of P. Neal \cite{N} (in the spirit of \cite{B-H-J})
in order to obtain the appropriate normalization of the random variable $T_{m}(N)$ and arrive into
specific limiting distributions. Our main results are Theorem 4.2 (see also Subsection 4.3, whose content complements Theorem 4.2). It is notable that for the considered class of coupon
probabilities the random variable $T_{m}(N)$, appropriate normalized, converges in distribution to a Gumbel random variable. This is a generalization of the classical result of P. Erd\H{o}s and A. R$\acute{e}$nyi (see (\ref{3})) for the case of equal coupon probabilities. For the special case $m=1$ the statement of Theorem 4.2 had been established in our earlier work \cite{DP}. Subsection 4.3 is a brief discussion of the case of slowly decaying
sequences. As an illustration we consider the sequence $a_j = 1 / j^p$, where $p > 0$. Such sequences were never studied before, not even in the case $m=1$. An interesting phenomenon arises regarding decaying sequences: If the decay of $\alpha$ is subexponential, then, at least for the great variety of cases we have analyzed, the limit distribution of $T_{m}(N)$ is always Gumbel.
However, if $\alpha$ decays to zero exponentially, then the limit distribution of $T_{m}(N)$ is not Gumbel and we expect that the same is true for
sequences decayng to zero superexponentially. In the latter case the behavior of $T_{m}(N)$ seems similar to the corresponding behavior for the case where $\alpha$ tends to infinity.

The proofs of Theorem 2.1 and of some technical lemmas of Section 2 are given in the Appendix (Section 5).
Finally, Section 6 is a short epilogue containing some concluding remarks and a comparison with earlier works.

\subsection{Some conjectures}

We finish this introductory section with two conjectures. Formulas \eqref{3}--\eqref{333} suggest the following conjecture:

\smallskip

\textbf{Conjecture 1.} In the case of equal coupon probabilities we have
\begin{equation}
V\left[\,T_m(N)\,\right] \sim \frac{\pi^2}{6} \, N^2,
\qquad
N \to \infty.
\label{SE22c}
\end{equation}
Actually, for $m=1$, we have already seen formula \eqref{4aa}, proved in \cite{B}, which is a stronger form of \eqref{SE22c}.

%
%In support of Conjecture 2 we present the following proposition:
%
%\smallskip
%
%\textbf{Proposition 1.} For fixed positive integers $m$ and $N$, the case of equal probabilities has the property that, among all choices of
%coupon probabilities, it is the one with the smallest $r$-th rising moment of  $T_{m}\left(N\right)$, for each $r$, as well as the one with the
%greatest moment generating function $E[z^{-T_m(N)}]$, for any $z > 1$.
%
%\smallskip
%
%\textit{Proof.} In view of (\ref{P7}), it is enough to show that for any fixed $t > 0$, the maximum of the quantity
%\begin{equation*}
%\prod_{j=1}^{N}\left[1-e^{-p_{j}t}\,S_{m}\left(p_{j}t \right) \right],
%\end{equation*}
%subject to the constraint
%\begin{equation}
%g(p_1, \dots, g_N) := p_{1}+\cdot \cdot \cdot +p_{N}=1,\qquad p_{j} \leq 0, \,\,j=1,2,\cdots,N,
%\label{A1}
%\end{equation}
%occurs when all $p_{j}$'s are equal. This can be done easily using Lagrange multipliers.
%\hfill $\blacksquare$
%
\smallskip

\textbf{Conjecture 2.} For fixed positive integers $m$ and $N$, the case of equal probabilities, has the property that it is the one with the smallest variance
$V\left[\, T_m(N)\,\right]$.

\smallskip

The results of the present work confirm that, for a large class of probabilities, $V\left[ T_{m}(N)\right]$ is actually minimized in the case of equal probabilities, as $N$ becomes sufficiently large.
% In particular, the variance $V\left[\, T_m(N)\,\right]$ is a Schur-convex function.

\section{Asymptotics of $E[T_m(N)]$ and $V[T_m(N)]$}

\subsection{Preliminary material}

When $N$ is large it is not clear at all what information one can obtain
for $E[T_m(N)]$ and $V[T_m(N)]$ from the formulas (\ref{5}) and (\ref{5A}) respectively. For this reason there is a need to develop efficient
ways for deriving asymptotics as $N \rightarrow \infty $. As in \cite{BP}, \cite{DP}, and \cite{DPM}, let $\alpha =\{a_{j}\}_{j=1}^{\infty }$ be a sequence of strictly positive
numbers. Then, for each integer $N > 0$, one can create a probability measure
$\pi _N =\{p_1,...,p_N\}$ on the set of types $\{1,...,N\}$ by taking
\begin{equation}
p_j = \frac{a_j}{A_N},
\qquad \text{where}\quad
A_N = \sum_{j=1}^N a_j.
\label{8}
\end{equation}
Notice that $p_j$ depends on $\alpha $ and $N$, thus, given $\alpha $, it
makes sense to consider the asymptotic behavior of $E\left[\, T_m(N)\,\right]$, $E\left[\,T_m(N)\left(T_{m}(N)+1\right)\,\right]$, and $V\left[\, T_m(N)\,\right]$ as $N\rightarrow \infty$.

\smallskip

\textbf{Remark 1.} Clearly, for given $N$ the $p_j$'s can be assumed monotone in $j$ without loss of generality. This tells us that if
$\alpha = \{a_j\}_{j=1}^{\infty }$ is such that $\lim_j a_j = \infty$, by rearranging its terms can be assumed without loss of
generality that $\alpha$ is a nondecreasing sequence. Likewise, if $\alpha = \{a_j\}_{j=1}^{\infty }$ is such that $\lim_j a_j = 0$,
then by rearranging its terms can be assumed without loss of generality that $\alpha$ is a nonincreasing sequence.

\smallskip

Motivated by (\ref{5}) we introduce the notation
\begin{align}
E_{m}(N;\alpha ):&=\int_{0}^{\infty}\left[1-\prod_{j=1}^{N}\bigg(1-e^{-a_{j}t}\,S_{m}\left(a_{j}t\right)\bigg)\right]dt \label{9} \\
&=\int_{0}^{1}\left[1-\prod_{j=1}^{N}\bigg(1-x^{a_{j}}\,S_{m}\left(-a_{j}\ln x\right)\bigg)\right]\frac{dx}{x}.
\label{10}
\end{align}
For a sequence $\alpha =\{a_{j}\}_{j=1}^{\infty }$ and a number $s > 0$ we set $s \alpha = \{sa_{j}\}_{j=1}^{\infty }$
(notice that $\alpha$ and $s \alpha$ create the same sequence of probability measures $\pi_N$, $N = 2, 3, \dots$).
Then, (\ref{9}) implies that
$E_{m}(N;s\alpha) = s^{-1} E_{m}(N;\alpha )$
and hence, in view of (\ref{5}) and (\ref{8}),
\begin{equation}
E\left[\, T_m(N)\,\right] =A_{N}\,E_{m}(N;\alpha).
\label{12}
\end{equation}
Likewise, motivated by (\ref{5A}), let us introduce
\begin{align}
Q_{m}(N;\alpha ):&=2\int_{0}^{\infty}t\left[1-\prod_{j=1}^{N}\bigg(1-e^{-a_{j}t}\,S_{m}\left(a_{j}t\right)\bigg)\right]dt \label{13}\\
&=-2\int_{0}^{1}\left[1-\prod_{j=1}^{N}\bigg(1-x^{a_{j}}\,S_{m}\left(-a_{j}\ln x\right)\bigg)\right]\frac{\ln x}{x}dx.\label{14}
\end{align}
From the above it follows that $Q_{m}(N;s\alpha) = s^{-2} Q_{m}(N;\alpha ),$ hence
\begin{equation}
E\left[\,T_m(N)\left(T_{m}(N)+1\right)\,\right] =A^{2}_{N}\,Q_{m}(N;\alpha).
\label{15}
\end{equation}
In view of (\ref{12}) and (\ref{15}), (\ref{VAR}) yields
\begin{equation}
V\left[\, T_m(N)\,\right]= A^{2}_{N}Q_{m}(N;\alpha)-A_{N}E_{m}(N;\alpha)-A^{2}_{N}E_{m}(N;\alpha)^{2}
\label{17a}.
\end{equation}
Under (\ref{8}) the problem of estimating $E\left[\, T_m(N)\,\right]$ can be treated as two separate problems, namely estimating $A_{N}$ and estimating $E_{m}(N;\alpha)$, (see (\ref{12})). The estimation of $A_N$ can be considered an
external matter which can be handled by existing powerful methods, such as the Euler-Maclaurin sum formula, the Laplace method for sums (see, e.g., \cite{B-O}), or even summation by parts. Hence, our analysis focuses on estimating $E_{m}(N;\alpha)$. Of course, the same observation applies to the
expression of (\ref{15}).

\subsection{The Dichotomy}

For convenience, we set
\begin{equation*}
f_{N}^{\alpha }(x) := \prod_{j=1}^{N}\bigg[1-x^{a_{j}}\,S_{m}\left(-a_{j}\ln x\right)\bigg],
\qquad 0 \leq x \leq 1,\\
\end{equation*}
in particular, $f_N^\alpha (0) := f_N^\alpha (0+) = 1$ and $ f_N^\alpha (1) = 0$.

Since
\begin{equation}
\frac{d}{dy} \left[e^{-y} S_m(y)\right] = -\frac{y^{m-1} e^{-y}}{(m-1)!},
\label{16aa}
\end{equation}
we get that $f_{N}^{\alpha }(x) $ is monotone decreasing in $x$. Furthermore, (\ref{7aa}) implies immediately that
$f_{N+1}^{\alpha }(x)\leq f_{N}^{\alpha }(x)$. In particular
\begin{equation*}
\lim_{N}f_{N}^{\alpha }(x)=\prod_{j=1}^{\infty}\bigg[1-x^{a_{j}}\,S_{m}\left(-a_{j}\ln x\right)\bigg] \qquad\text{exists for all }\; x \in [0, 1].
\end{equation*}
Thus, by applying the Monotone Convergence Theorem in (\ref{10}) and (\ref{14}), we get respectively
\begin{equation}
L_{1}(\alpha;m ):=\lim_{N}E_{m}(N;\alpha )=\int_{0}^{1}\left\{1-\prod_{j=1}^{\infty}\bigg[1-x^{a_{j}}\,S_{m}(-a_{j}\ln x)\bigg]\right\}
\frac{dx}{x}
\label{17}
\end{equation}
and
\begin{equation}
L_{2}(\alpha;m ):=\lim_{N}Q_{m}(N;\alpha )=-2\int_{0}^{1}\left\{1-\prod_{j=1}^{\infty}\bigg[1-x^{a_{j}}\,S_{m}\left(-a_{j}\ln x\right)\bigg]\right\}\frac{\ln x}{x}\, dx.
\label{18}
\end{equation}
Notice that $L_{1}(\alpha ;m ),L_{2}(\alpha;m ) > 0$, for any $\alpha $ (since, for
every $x\in (0,1)$, $f_{N}^{\alpha }(x)<1$ and decreases with $N$).
However, we may have $L_{1}(\alpha;m )=\infty $ and/or $L_{2}(\alpha;m )=\infty $. In fact, Theorem 2.1 below tells us that
$L_{1}(\alpha;m )= \infty$ if and only if $L_{2}(\alpha;m ) = \infty$.

\smallskip
\begin{theorem}
Let $L_1 (\alpha;m)$ and $L_2 (\alpha;m)$ as defined in (\ref{17}) and (\ref{18}) respectively. The following are equivalent
(for all positive integers $m$):\\
(i) $L_1 (\alpha;m) < \infty$ \\
(ii) $L_2 (\alpha;m) < \infty$\\
(iii) There exist a $\xi \in (0,1)$ such that
\begin{equation}
\sum_{j = 1}^\infty \xi^{a_j} < \infty.
\label{heat}
\end{equation}
\end{theorem}
The proof of the theorem is given in the Appendix.

The theorem implies that we have the following dichotomy simultaneously for all positive integers $m$:
\begin{equation}
\text{(i) }\; 0<L_{1}(\alpha;m ) ,\; L_{2}(\alpha;m )<\infty
\quad \text{or \quad
(ii) }\; L_{1}(\alpha;m )=L_{2}(\alpha;m )=\infty.
\label{B2}
\end{equation}

\textbf{Remark 2.} The word ``dichotomy" may be misleading: If $p > 0$, then, the sequence $\alpha = \{e^{pj} \}_{j=1}^{\infty}$
satisfies $L_i(\alpha; m) < \infty$, while for the sequence $\beta = \{e^{-pj} \}_{j=1}^{\infty}$ we
have $L_i(\beta; m) = \infty$ ($i = 1, 2$). However, it is clear that $\alpha$ and $\beta$ produce the same coupon probabilities(!),
i.e. the same sequence of probability measures $\{\pi_N\}_{N=2}^{\infty}$. This is an exceptional case, since apart from this case, it can
be shown by straightforward induction on $N$ that if
$\alpha = \{a_j \}_{j=1}^{\infty}$ is a sequence such that $\lim_j a_j = \infty$, then there is no sequence $\beta = \{b_j \}_{j=1}^{\infty}$,
with $\lim_j b_j = 0$, producing the same coupon probabilities as $\alpha$ and vice versa.

\subsection{Case I: $L_1(\protect\alpha; m) < \infty$}

Let $A_{N}$ and $L_1(\alpha;m )$ be as in (\ref{8}) and (\ref{17}). We note that, by
Theorem 2.1 (see (\ref{heat})), $L_1(\alpha; m) < \infty$ implies that $\lim_{j} a_j = \infty$ (hence $\lim_N A_N = \infty$).

\smallskip

\begin{theorem}
If $L_1(\alpha; m) < \infty $, then for all positive integers $m$ we have
\begin {align}
E\left[\,T_{m}(N)\,\right] &=A_{N}L_{1}(\alpha;m )\left( 1 + \delta_N\right),
\label{fm}
\\
E\left[\,T_{m}(N)(T_{m}(N)+1)\,\right] &=A_{N}^{2}L_{2}(\alpha;m )\left( 1 + \Delta _N\right),
\label{sm}
\\
V\left[\, T_m(N)\,\right]&=A_{N}^{2}\left[L_{2}(\alpha;m)-L_{1}(\alpha;m)^{2}\right]\left[1 + o\left(1\right)\right]
\label{Bb}
\end{align}
(as $N \rightarrow \infty$), where for the error terms
\begin{equation}
\delta_N := L_{1}(\alpha;m ) - E_{m}(N;\alpha)
\quad\text{and}\quad
\Delta _N := L_{2}(\alpha;m ) - Q_{m}(N;\alpha)
\label{ET1}
\end{equation}
we have $\delta_N = o(1)$ and $\Delta_N = o(1)$ as $N \rightarrow \infty$.
Furthermore, in the formula \eqref{Bb} it is always true that
\begin{equation}
L_2(\alpha; m) - L_1(\alpha; m)^2 > 0.
\label{8bb}
\end{equation}
\end{theorem}

{\it Proof}. Formula (\ref{fm}) follows immediately from (\ref{12}) and (\ref{17}), while formula (\ref{sm})
follows from (\ref{15}) and (\ref{18}).

To prove (\ref{8bb}) we first notice that
\begin{equation*}
G(x) := 1 - \prod_{j=1}^{\infty} \bigg[1-x^{a_{j}}\,S_{m}\left(-a_{j}\ln x\right)\bigg]
\end{equation*}
is a nondegenerate distribution function on $[0, 1]$. If $X$ is a random variable with distribution function $G(x)$, then simple
integration by parts in (\ref{18}) and (\ref{17}) gives
\begin{equation*}
L_2(\alpha; m) = E_G\left[\ln(x)^2\right] > E_G\left[\ln(x)\right]^2 = L_1(\alpha; m)^2,
\end{equation*}
where $E_G[\,\cdot\,]$ denotes the expectation associated to the distribution function $G(x)$. Having established (\ref{8bb}), formula (\ref{Bb})
follows by using (\ref{fm}) and (\ref{sm}) in (\ref{VAR}).
\hfill $\blacksquare$

\smallskip

\textbf{Remark 3.} If $a_j$ grows to infinity sufficiently fast, we can get a better
estimate for the errors $\delta_N$ and $\Delta_N$ of \eqref{ET1}.
By (\ref{10}), (\ref{17}), (\ref{LL1b}) (see Appendix), and Tonelli's Theorem
\begin{align*}
\delta _{N}&=\int_{0}^{1}\prod_{j=1}^{N}\bigg(1-x^{a_{j}}S_{m}\left(-a_{j}\ln x\right)\bigg)\left[1-\prod_{j=N+1}^{\infty}\bigg(1-x^{a_{j}}S_{m}\left(-a_{j}\ln x\right)\bigg)\right]\frac{dx}{x}\nonumber\\
&\leq \sum_{j=N+1}^{\infty}\left[\int_{0}^{1}x^{a_{j}}S_{m}\left(-a_{j}\ln x\right)\frac{dx}{x}\right]=\sum_{j=N+1}^{\infty}\sum_{k=0}^{m-1}\int_{0}^{1}x^{a_{j}-1}\left(\ln x\right)^{k}.
\end{align*}
Integration by parts yields
\begin{equation}
\delta _{N}\leq m\sum_{j=N+1}^{\infty}\frac{1}{a_{j}}.
\label{ET2}
\end{equation}
In a similar manner one gets
\begin{equation}
\Delta _{N}\leq m\left(m+1\right)\sum_{j=N+1}^{\infty}\frac{1}{a^{2}_{j}}.
\label{ET3}
\end{equation}

\smallskip

\textbf{Remark 4.} For $r = 1, 2, \dots$ let us set
\begin{equation*}
L_r(\alpha; m) := (-1)^{r-1} r\int_{0}^{1}\left\{1-\prod_{j=1}^{\infty}\bigg[1-x^{a_{j}}\,S_{m}\left(-a_{j}\ln x\right)\bigg]\right\}\frac{\ln^{r-1} x}{x} \,dx,
%\label{18}
\end{equation*}
i.e.
\begin{equation*}
L_r(\alpha; m) = r\int_0^{\infty} \left\{1-\prod_{j=1}^{\infty}\left[1 - S_{m}(a_j t) e^{-a_j t} \right]\right\} t^{r-1} dt.
%\label{18}
\end{equation*}
Then, Theorem 2.1 is valid for $L_r(\alpha; m)$, for any $r$ (the proof is similar). Furthermore, it is not hard to see that
\begin{equation*}
E\left[ T_m(N)^{(r)} \right] = A_N^{r} L_{r}(\alpha; m)\left[1 + o(1)\right],
\qquad
N \to \infty,
\end{equation*}
which is an extension of Theorem 2.2 for all $r$.

\subsection{Case II: $L_1(\protect\alpha; m) = \infty$}

\subsubsection{Asymptotic Behavior of $E[T_m(N)]$}

By Theorem 2.1, $L_1(\alpha;m )=\infty $ is equivalent to
\begin{equation*}
\sum_{j=1}^{\infty }x^{a_{j}}=\infty \text{,}
\qquad
\text{for all } x \in (0,1).
\end{equation*}
For our further analysis we follow \cite{BP}, \cite{DP}, and \cite{DPM}, and write $a_{j}$ in the form
\begin{equation}
a_{j}=\frac{1}{f(j)},
\label{aj}
\end{equation}
where
\begin{equation}
f(x) > 0 \qquad \text{and} \qquad f'(x) > 0.
\label{C1a}
\end{equation}
In order to proceed we assume that $f(x)$ possesses three derivatives satisfying the following conditions as $x\rightarrow \infty$:
\begin{align}
\text{(i) }& f(x)\rightarrow \infty,
& &\text{(ii) } \frac{f^{\prime }(x)}{f(x)}\rightarrow 0,\nonumber \\
\text{(iii) } &\frac{f^{\prime \prime}(x)/f^{\prime }(x)}{f'(x)/f(x)} = O\left(1\right),
& &\text{(iv) } \frac{f^{\prime \prime\prime}(x)\;f(x)^{2}}{ f^{\prime }(x)^{3}} = O\left(1\right).
\label{C1}
\end{align}
Roughly speaking, $f(\cdot)$ belongs to the class of positive and strictly increasing $C^3(0, \infty)$ functions, which grow to $\infty$
(as $x \rightarrow \infty$) slower than exponentials, but faster than powers of logarithms.
These conditions are satisfied by a variety of commonly used functions. For example,
\begin{equation*}
f(x) = x^p (\ln x)^q, \quad p > 0,\ q \in \mathbb{R},\qquad \qquad
f(x) = \exp(x^{r}),\quad 0 < r < 1,
\end{equation*}
or various convex  combinations of products of such functions. Notice that the smoothness assumption on $f$ does not impose any restriction on
the sequence $\alpha$, since we only need $f(x)$ to interpolate $1 / a_j$ for $x = j$, $j = 1, 2, \dots \;$. The restrictions of $\alpha$ come from
the growth assumptions \eqref{C1}. In Subsection 4.3 we discuss the case where $f(x)$ grows slower than any (positive) power of $x$ and hence does not
satisfy all conditions of \eqref{C1}.

For typographical convenience we set
\begin{equation}
F(x) := f(x) \ln\left(\frac{f(x)}{f'(x)}\right)
\label{F}
\end{equation}
(notice that (\ref{C1a}) and (ii) of (\ref{C1}) imply that $F(x) > 0$, for $x$ sufficiently large). The following lemma plays an important role in our analysis:

\smallskip

\begin{lemma}
Set
\begin{equation}
J_{\kappa}(N):=\int_{1}^{N} f(x)^{\kappa}e^{-\frac{F(N)}{f(x)}s}dx ,\quad \kappa \in \mathbb{R}.
\label{I}
\end{equation}
Then, under (\ref{C1}) and (\ref{F}), we have
\begin{equation}
J_{\kappa}(N)=\frac{f(N)^{\kappa+2}}{sF(N)f^{\prime}(N)}e^{-\frac{F(N)}{f(N)}s}
+\omega(N)\;\frac{f(N)^{\kappa+3}}{s^{2}F(N)^{2}f^{\prime}(N)}e^{-\frac{F(N)}{f(N)}s}\left[1+O\left(\frac{f(N)}{F(N)}\right)\right], \label{IIa}
\end{equation}
where
\begin{equation}
\omega(N) := -2+\frac{f^{\prime \prime}(N)/f^{\prime }(N)}{f^{\prime }(N)/f(N)},
\label{a}
\end{equation}
uniformly in $s\in [s_{0},\infty)$, for any fixed $s_{0} > 0 $.
\end{lemma}
The proof is given in \cite{DP} in the case where $\kappa \geq 0$, while
it is straightforward to check that the lemma above is still valid when $\kappa$ is negative. Notice that the condition (iii) of (\ref{C1}) says that
$\omega(N) = O(1)$ as $N \to \infty$.\\
Using Lemma 2.1, as well as (\ref{C1a}), (\ref{C1}), and (\ref{F}), we have as $N\rightarrow \infty$
\begin{equation}
\int_{1}^{N}e^{-\frac{F(N)}{f(x)}s}S_{m}\left(\frac{F(N)}{f(x)}s\right)dx
\sim \frac{1}{\left(m-1\right)!}
\left[\ln \left(\frac{f(N)}{f^{\prime }(N)}\right)\right]^{m-2}s^{m-2}\left[ \frac{f(N)}{f^{\prime }(N)}\right] ^{1-s},
\label{C5}
\end{equation}
where, as usual, $E_1(N) \sim E_2(N)$ means that $E_1(N) / E_2(N) \to \infty$ as $N\rightarrow \infty$.
From (\ref{C5}) we obtain
\begin{equation}
\lim_{N}\int_{1}^{N}e^{-\frac{F(N)}{f(x)}s}S_{1}\left(\frac{F(N)}{f(x)}s\right)dx =\left\{
\begin{array}{rll}
\infty,&  \text{if }& s<1, \\
0,&   \text{if }& s\geq 1,
\end{array}
\right.  \label{LL2a}
\end{equation}
\begin{equation}
\lim_{N}\int_{1}^{N}e^{-\frac{F(N)}{f(x)}s}S_{2}\left(\frac{F(N)}{f(x)}s\right)dx =\left\{
\begin{array}{rcc}
\infty,&  \text{if } s<1, \\
1,& \text{if } s=1, \\
0,&   \text{if } s> 1,
\end{array}
\right.
\label{LL1}
\end{equation}
while, for $m\geq3$
\begin{equation}
\lim_{N}\int_{1}^{N}e^{-\frac{F(N)}{f(x)}s}S_{m}\left(\frac{F(N)}{f(x)}s\right)dx =\left\{
\begin{array}{rll}
\infty,&  \text{if }& s\leq1, \\
0,&   \text{if }& s> 1.
\end{array}
\right.  \label{LL2}
\end{equation}
It is easy for one to check that the function $h(x):=e^{-\frac{F(N)}{f(x)}s}S_{m}\left(\frac{F(N)}{f(x)}s\right)$ is increasing. Hence,
\begin{align*}
\int_{1}^{N}e^{-\frac{F(N)}{f(x)}s}S_{m}\left(\frac{F(N)}{f(x)}s\right)dx&\leq \sum_{j=1}^{N}e^{-\frac{F(N)}{f(j)}s}S_{m}\left(\frac{F(N)}{f(j)}s\right) \nonumber\\
&\leq \int_{1}^{N}e^{-\frac{F(N)}{f(x)}s}S_{m}\left(\frac{F(N)}{f(x)}s\right)dx \nonumber\\
&\,\,\,\,\,\,+e^{-\frac{F(N)}{f(N)}s}S_{m}\left(\frac{F(N)}{f\left(N\right)}s\right).
\end{align*}
It follows (see (\ref{F}) and (ii) of (\ref{C1})) that the limits in (\ref{LL1}) and (\ref{LL2}) are valid, if the integral is replaced by the sum, namely
$\sum_{j=1}^{N}e^{-\frac{F(N)}{f(j)}s}S_{m}\left(\frac{F(N)}{f(j)}s\right)$.
Finally, by the definition of $F(\cdot)$ and the Taylor expansion for the logarithm, namely $\ln(1-x)\sim-x$ as $x\rightarrow0$, we get
\begin{equation}
\lim_{N}\sum_{j=1}^{N}\ln \Bigg( 1-e^{-\frac{F(N)}{f(j)}s}S_{1}\left(\frac{F(N)}{f(j)}s\right)\Bigg) =\left\{
\begin{array}{rc}
-\infty,&  \text{if } s<1 \\
0,&   \text{if } s\geq1,
\end{array}
\right.  \label{SL1A}
\end{equation}
\begin{equation}
\lim_{N}\sum_{j=1}^{N}\ln \Bigg( 1-e^{-\frac{F(N)}{f(j)}s}S_{2}\left(\frac{F(N)}{f(j)}s\right)\Bigg) =\left\{
\begin{array}{rc}
-\infty,&  \text{if } s<1 \\
0,&   \text{if } s>1\\
-1,&\text{if } s=1,
\end{array}
\right.  \label{SL1}
\end{equation}
\begin{equation}
\lim_{N}\sum_{j=1}^{N}\ln \Bigg( 1-e^{-\frac{F(N)}{f(j)}s}S_{m}\left(\frac{F(N)}{f(j)}s\right)\Bigg) =\left\{\begin{array}{rc}
-\infty,&  \text{if } s\leq1 \\
0,&   \text{if } s> 1,
\end{array}
\right.  \label{SL2}
\end{equation}
for $m=1$, $m=2$, and $m=3, 4, \cdots$ respectively. Next, we take advantage of the above limits. Starting from (\ref{9}), and for any given $\varepsilon \in (0,1)$ we rewrite $E_{m}(N;\alpha)$ as
\begin{align}
E_{m}(N;\alpha)= F(N)\left[\,1+\varepsilon -I_1 (N)-I_2 (N)+I_3 (N)\,\right],
\label{b3a}
\end{align}
where
\begin{align}
I_1 (N):&= \int_0^{1-\varepsilon} \left[\exp \left\{ \sum_{j=1}^{N}\ln \Bigg( 1-e^{-\frac{F(N)}{f(j)}s}S_{m}\left(\frac{F(N)}{f(j)}s\right)\Bigg) \right\}\right] ds,\label{I1}\\
I_{2}(N): &= \int_{1-\varepsilon}^{1+\varepsilon }\left[ \exp \left\{ \sum_{j=1}^N \ln \Bigg( 1-e^{-\frac{F(N)}{f(j)}s}S_{m}\left(\frac{F(N)}{f(j)}s\right)\Bigg) \right\} \right] ds,\label{I2}\\
I_{3}(N): &=\int_{1+\varepsilon}^{\infty }\left[ 1-\exp \left\{ \sum_{j=1}^N \ln \Bigg( 1-e^{-\frac{F(N)}{f(j)}s}S_{m}\left(\frac{F(N)}{f(j)}s\right)\Bigg) \right\} \right] ds.\label{I3}
\end{align}
For typographical convenience we set
\begin{equation}
\delta := \frac{1}{\ln \left(\frac{f(N)}{f^{\prime}(N)}\right)}=\frac{f(N)}{F(N)}
\label{delta}
\end{equation}
(notice that as $N \rightarrow \infty$, $\delta \rightarrow 0^+$).

\smallskip

\begin{lemma}
Let $I_{1}(N)$, $I_{2}(N)$, $I_{3}(N)$, and $\delta$ be as defined in (\ref{I1}), (\ref{I2}), (\ref{I3}), and (\ref{delta})
respectively. Then, for any given $\varepsilon \in (0,1)$ and for all positive integers $m$ we have, as $\delta \rightarrow 0^{+}$,
\begin{equation}
I_{1}(N) = o\left( \delta^{4-m}e^{-\varepsilon/ \delta}\right).
\label{b3}
\end{equation}
Furthermore,
\begin{align}
I_{2}(N)=&\,\varepsilon+\left(m-2\right)\delta \ln \delta+\left[\,\ln\left(m-1\right)!- \gamma\,\right] \delta
+\left(m-2\right)^{2} \delta^{2}\ln \delta\nonumber\\
&\,\,\,\,+\left[\left(m-2\right)\ln\left(m-1\right)!
-\left(m-2\right)\gamma-\left(m-1\right)-\omega(N)\left(m-1\right)!\right]\delta^{2}\nonumber\\
&\,\,\,\,+O\left(\delta^{3}\left(\ln \delta\right)^{2}\right).\label{eint11}
\end{align}
and
\begin{equation}
I_{3}(N)=\frac{\left(1+\varepsilon\right)^{m-2}}{\left(m-1\right)!}\frac{1}{\delta^{m-3}}\,e^{-\varepsilon/\delta}\left(1+O\left(\delta\right)\right) \quad \text{as}\,\,\, \delta\rightarrow 0^{+}. \label{bb3}
\end{equation}
\end{lemma}
The proof of the lemma is given in the Appendix.

\smallskip

\textbf{Observation 1.} It follows by Lemma 2.2 that both integrals $I_{1}(N)$ and $I_{3}(N)$, are negligible compared to the sixth term in the asymptotic expansion of the integral $I_{2}(N)$. Hence, all the information for the $E\left[\, T_m(N)\,\right]$, at least for the five first terms, comes from $I_{2}(N)$.

\smallskip

We are, therefore, ready to present the following theorem.

\smallskip

\begin{theorem} Let $\delta$ be as defined in (\ref{delta}) (hence $\delta\rightarrow 0^{+}$ as $N\rightarrow \infty$) and  $\omega(N)$
as given in (\ref{a}). Then ($\gamma$ is, as usual, the Euler-Mascheroni constant)
\begin{align}
&E\left[\, T_m(N)\,\right] = A_N f(N) \left\{\frac{1}{\delta}-\left(m-2\right)\ln \delta + [\,\gamma-\ln\left(m-1\right)!\,]
-(m-2)^2 \delta\ln \delta\right.
\nonumber
\\
& +[(m-1) + \omega(N) (m-1)! - (m-2) \ln (m-1)! +(m-2) \gamma]\delta
+ O\left(\delta^2 \ln^2\delta \right) \bigg\}.
\label{R1}
\end{align}
\end{theorem}
{\it Proof}. The result follows immediately by combining (\ref{12}), (\ref{b3a}), and Lemma 2.2.
\hfill $\blacksquare$

\smallskip

To follow D.J. Newman and L. Shepp \cite{N-S}, although the first set ``costs'' $A_{N}f(N)/\delta$, all further sets cost $A_{N}f(N)\ln\delta$.

\subsubsection{Asymptotics of the second rising moment of $T_m(N)$}

We will follow a similar approach as in Subsubsection 2.4.1, in order to find the sixth(!) term in the asymptotic expansion of the second rising moment of the random variable $T_m(N)$, so that the leading behavior of $V[\,T_{m}(N)\,]$ will be obtained. Let us expand $Q_m(N;\alpha)$ as
\begin{align}
Q_{m}(N;\alpha)= 2F(N)^{2}\left[\frac{1}{2}+\varepsilon+\frac{\varepsilon^{2}}{2}-I_{4}(N)-I_{5}(N) + I_{6}(N)\right],\label{QM}
\end{align}
where
\begin{equation}
I_4 (N):= \int_0^{1-\varepsilon} \left[\exp \left\{ \sum_{j=1}^{N}\ln \Bigg( 1-e^{-\frac{F(N)}{f(j)}s}S_{m}\left(\frac{F(N)}{f(j)}s\right)\Bigg) \right\}\right]s\, ds,
\label{I4}
\end{equation}
\begin{equation}
I_{5}(N) = \int_{1-\varepsilon}^{1+\varepsilon }\left[ \exp \left\{ \sum_{j=1}^N \ln \Bigg( 1-e^{-\frac{F(N)}{f(j)}s}S_{m}\left(\frac{F(N)}{f(j)}s\right)\Bigg) \right\} \right] s\,ds,
\label{I5}
\end{equation}
and
\begin{equation}
I_{6}(N) =\int_{1+\varepsilon}^{\infty }\left[ 1-\exp \left\{ \sum_{j=1}^N \ln \Bigg( 1-e^{-\frac{F(N)}{f(j)}s}S_{m}\left(\frac{F(N)}{f(j)}s\right)\Bigg) \right\} \right] s\,ds.
\label{I6}
\end{equation}

\begin{lemma}
Let $I_{4}(N)$, $I_{5}(N)$, $I_{6}(N)$, and $\delta$ be as defined in (\ref{I4}), (\ref{I5}), (\ref{I6}), and (\ref{delta}) respectively. Then, for any given $\varepsilon \in (0,1)$ we have, as $\delta\rightarrow 0^{+}$,
\begin{equation}
I_{4}(N) = o\left(\delta^{4-m}e^{-\varepsilon/ \delta}\right).
\label{b32}
\end{equation}
Furthermore,
\begin{align}
I_{5}(N)=&\; \varepsilon + \frac{\varepsilon^{2}}{2}+\left(m-2\right)\delta\ln \delta+\left[\ln\left(m-1\right)!-\gamma\right]\delta
-\frac{\left(m-2\right)^{2}}{2} \delta^{2}\ln^{2} \delta\nonumber\\
&\,\,\,+\left[\left(m-2\right)^{2}-\left(m-2\right)\left(\ln\left(m-1\right)!-\gamma\right)\right]\delta^{2}\ln \delta\nonumber\\
&\,\,\,+\left[\left(m-2\right)\ln\left(m-1\right)!-\left(m-2\right)\gamma-\omega(N)\left(m-1\right)!-\left(m-1\right)\right.\nonumber\\
&\left.\,\,\,\,\,\,\,\,\,-\frac{1}{2}\left(\ln\left(m-1\right)!\right)^{2}-\frac{1}{2}\left(\gamma^{2}+\frac{\pi^{2}}{6}\right)+\gamma \ln\left(m-1\right)!\right]\delta^{2}+O\left(\delta^{3}\left(\ln\delta\right)^{2}\right)
\label{Beint22}
\end{align}
and
\begin{equation}
I_{6}(N)=\frac{\left(1+\varepsilon\right)^{m-1}}{\left(m-1\right)!}\frac{1}{\delta^{m-3}}\,e^{-\varepsilon/\delta}\left(1+O\left(\delta\right)\right). \label{bb32}
\end{equation}
\end{lemma}
The proof of lemma above is given in the Appendix.

\smallskip

\textbf{Observation 2.} By Lemma 2.3 we have that both integrals $I_{4}(N)$ and $I_{6}(N)$, are negligible compared to the seventh term in the asymptotic expansion of the integral $I_{5}(N)$. Hence, all the information regarding $E\left[\,T_m(N)\left(T_{m}(N)+1\right)\,\right]$, at least for the six first terms, comes from $I_{5}(N)$.

\smallskip

By combining (\ref{15}), (\ref{QM}), and Lemma 2.3 we obtain the following theorem.

\smallskip

\begin{theorem}
Let $\delta$ be as defined in (\ref{delta}) (hence $\delta\rightarrow 0^{+}$ as $N\rightarrow \infty$) and  $\omega(N)$ as given in (\ref{a}). Then for all positive integers $m$
\begin{align}
E\left[\,T_m(N)\left(T_{m}(N)+1\right)\,\right]& = A^2_N f(N)^2 \left\{\frac{1}{\delta^{2}}-2\left(m-2\right)\frac{\ln\delta}{\delta}
-2\left[\ln\left(m-1\right)!-\gamma\right]\frac{1}{\delta}\right.
\nonumber
\\
& + (m-2)^2 \ln^2\delta + 2(m-2) [\ln(m-1)!-\gamma - (m-2)] \ln\delta
\nonumber
\\
& +[2\left(m-2\right)\gamma-2\left(m-2\right)\ln\left(m-1\right)!+2\,\omega(N)(m-1)!
\nonumber
\\
& +2\left(m-1\right)+\left(\ln\left(m-1\right)!\right)^{2}+\gamma^2 + (\pi^2 / 6) - 2\gamma \ln\left(m-1\right)!]
\nonumber\\
& +O\left(\delta \ln^2\delta \right) \bigg\}.
\label{R2}
\end{align}
\end{theorem}
We are now ready for our main result regarding the variance (in Case II).

\subsubsection{Asymptotics of $V\left[T_{m}(N)\right]$}

\begin{theorem} Let $\alpha = \{a_j\}_{j=1}^{\infty } = \{1/f(j)\}_{j=1}^{\infty}$, where $f$ satisfies (\ref{C1a}) and (\ref{C1})
(hence, $L_1(\protect\alpha;m ) = \infty$). Then for all positive integers $m$ we have as $N\rightarrow \infty$
\begin{equation}
V\left[\,T_m(N)\,\right] \sim \frac{\pi^2}{6}\;A^2_N\;f(N)^2
= \frac{\pi^2}{6} \cdot \frac{1}{p_N^2} = \frac{\pi^2}{6} \cdot \frac{1}{\min_{1 \leq j \leq N}\, \{p_j\}^2},
\label{FINAL}
\end{equation}
where $A_N = \sum_{j=1}^N a_j$ ($p_j = a_j / A_N$ are the coupon probabilities).
\end{theorem}
{\it Proof}. From Theorems 2.3 and 2.4 one gets
\begin{equation*}
E\left[\,T_{m}(N)(T_{m}(N)+1)\,\right]-E\left[\,T_{m}(N)\,\right]^{2}\sim \frac{\pi^{2}}{6}\;A^{2}_{N}\;f(N)^{2}\;\;\;\text{as $N\rightarrow \infty$.}
\end{equation*}
In view of (\ref{17a}), in order to finish the proof it only remains to show that
\begin{equation}
\frac{E\left[\,T_m(N)\,\right]}{A^{2}_{N}\;f(N)^{2}}\rightarrow 0,\;\;\;N\rightarrow \infty.
\label{dom1}
\end{equation}
From (\ref{R1}) and  (\ref{delta})  we have
\begin{equation*}
E\left[\,T_m(N)\,\right] \sim A_{N}f(N)\ln\left(\frac{f(N)}{f^{\prime}(N)}\right).
\end{equation*}
Due to the above, (\ref{dom1}) is equivalent to
\begin{equation}
\frac{\ln f(N)-\ln f^{\prime}(N)}{A_{N}\;f(N)}\rightarrow 0,\;\;\;N\rightarrow \infty.
\label{dom3}
\end{equation}
Using (i) and (ii) of (\ref{C1}) it remains to prove that for sufficiently lage $x$
\begin{equation}
\frac{\ln f^{\prime}(x)}{\ln f(x)} = O(1).
\label{OOO}
\end{equation}
One arrives at (\ref{OOO}) starting from $\text{(iii)}$ (of (\ref{C1})). There is a positive constant $M,$ such that for sufficiently large $x$
\begin{equation*}
\left|\left(\ln f^{\prime}(x)\right)^{\prime}\right|\leq M\left|\left(\ln f(x)\right)^{\prime}\right|.
\end{equation*}
Since $f(x),\; f^{\prime}(x)>0$ the above becomes
\begin{equation*}
\left|\left(\ln f^{\prime}(x)\right)^{\prime}\right|\leq M\left(\ln f(x)\right)^{\prime}.
\end{equation*}
For any fixed $x_{0}>0$ and $x$ sufficiently large, we have
\begin{equation*}
\int_{x_{0}}^{x}\left|\left(\ln f^{\prime}(x)\right)^{\prime}\right|dx \leq M \int_{x_{0}}^{x}\left(\ln f(x)\right)^{\prime}.
\end{equation*}
Hence,
\begin{equation*}
\left|\int_{x_{0}}^{x}\left(\ln f^{\prime}(x)\right)^{\prime}dx\right| \leq M \int_{x_{0}}^{x}\left(\ln f(x)\right)^{\prime},
\end{equation*}
which implies
\begin{equation*}
\left|\ln f^{\prime}(x)-\ln f^{\prime}(x_{0})\right| \leq M \left(\ln f(x)-\ln f(x_{0})\right).
\end{equation*}
If we divide the above inequality with the positive function $\ln f(x)$ and use $\text{(i)}$ (of (\ref{C1}))  we have the desired result. This completes the proof.
\hfill $\blacksquare$

\textbf{Remark 5.} (a) It is notable that for the sequences of Case II the leading behavior of the variance of the random variable $T_{m}(N)$ is \textit{independent} of the value of the positive integer $m$ (which is in agreement with (\ref{4aa})). The reader may compare this with (\ref{Bb}) where the leading behavior of the variance \textit{depends} on $m$. \\
(b) Regarding the asymptotics of $A_N$ let us mention that if
\begin{equation*}
C_{f}:=\sum_{n=1}^{\infty}1/f(n) < \infty,
\end{equation*}
then
\begin{equation*}
A_{N}=C_{f}\left[1+o(1)\right].
\end{equation*}
On the other hand, if $C_{f}=\infty$, then, as $N\rightarrow \infty$, we have
\begin{equation*}
A_{N}\sim \int_{1}^{N}\frac{dx}{f(x)}.
\end{equation*}

\textbf{Remark 6.} Using the techniques presented in this section it can be shown that
\begin{equation*}
E\left[ T_m(N)^{(r)} \right] \sim A^r_{N} f(N)^{r} \ln \left(\frac {f(N)}{f^{\prime}(N)}\right)^r,
\qquad
r \in \mathbb{N},
\end{equation*}
where, as usual, the coupon probabilities are given by (\ref{8}) with
$\alpha = \{a_j\}_{j=1}^{\infty } = \{1/f(j)\}_{j=1}^{\infty}$, where $f$ satisfies (\ref{C1a}) and (\ref{C1}).

\section{Some examples}

\textbf{Example 1} (the positive power law)\textbf{.} Consider the sequence $\alpha = \{j^p\}_{j=1}^{\infty}$, where $p > 0$.
Here we have (see \eqref{17} and \eqref{18} and Theorem 2.1),
\begin{equation}
L_1(\alpha; m) = \int_0^{\infty} \left\{1 - \prod_{j=1}^{\infty}\bigg[1 - e^{-j^p t}\,S_{m}(j^p t)\bigg]\right\}dt < \infty
\label{SE5}
\end{equation}
and
\begin{equation}
L_2(\alpha; m) = 2 \int_0^{\infty} \left\{1 - \prod_{j=1}^{\infty}\bigg[1 - e^{-j^p t}\,S_{m}(j^p t)\bigg]\right\} t\, dt
< \infty,
\label{SE6}
\end{equation}
where $S_m(\cdot)$ is given by \eqref{7}. By Theorem 2.2 it follows that
\begin {align}
E\left[\,T_{m}(N)\,\right] &= A_N \, L_{1}(\alpha;m )\left( 1 + \delta_N \right),
\label{SE1}
\\
E\left[\,T_{m}(N)(T_{m}(N)+1)\,\right] &= A_N^2 \, L_{2}(\alpha;m )\left( 1 + \Delta_N\right),
\label{SE2}
\end{align}
where
\begin{equation}
\delta_N = o(1),
\qquad \qquad
\Delta_N = o(1)
\label{SE4a}
\end{equation}
and
\begin{equation}
A_N = \sum_{j=1}^N j^p = \frac{N^{p+1}}{p+1}\left[1 + O\left(\frac{1}{N}\right)\right],
\qquad
N \to \infty.
\label{SE4}
\end{equation}
Thus, \eqref{VAR} yields
\begin{equation}
V\left[\, T_m(N)\,\right] = \frac{N^{2p+2}}{(p+1)^2} \left[L_{2}(\alpha;m)-L_{1}(\alpha;m)^{2}\right]\left[1 + o\left(1\right)\right]
\label{SE3}
\end{equation}
as $N \to \infty$,
Actually, by \eqref{ET2} and \eqref{ET3} of Remark 3 we have
\begin{equation}
\delta _{N} \leq m\sum_{j=N+1}^{\infty}\frac{1}{j^p} = O\left(\frac{1}{N^{p-1}}\right),
\qquad \text{if }\; p > 1
\label{SE6}
\end{equation}
and
\begin{equation}
\Delta _{N}\leq m\left(m+1\right)\sum_{j=N+1}^{\infty}\frac{1}{j^{2p}} = O\left(\frac{1}{N^{2p-1}}\right),
\qquad \text{if }\; p > \frac{1}{2}
\label{SE7}
\end{equation}
(for the equalities in \eqref{SE6} and \eqref{SE7} see, e.g., \cite{A}). Thus, if $p > 1$, formula \eqref{VAR} gives
\begin{equation}
V\left[\, T_m(N)\,\right] = \frac{N^{2p+2}}{(p+1)^2} \left[L_{2}(\alpha;m)-L_{1}(\alpha;m)^{2}\right]\left[1 + O\left(\frac{1}{N^{(p-1)\wedge 1}}\right)\right],
\label{SE8}
\end{equation}
where $(p-1)\wedge 1 = \min\{p-1, 1\}$.

In the case $p = m = 1$ we can get explicit values for $L_1(\alpha; 1)$ and $L_2(\alpha; 1)$ as well as more accurate estimates for $\delta_N$
and $\Delta_N$ \cite{BP}, \cite{DP}.

\smallskip

\textbf{Example 2} (the generalized Zipf law)\textbf{.} Here, let us consider the sequence $\alpha = \{1 / j^p\}_{j=1}^{\infty}$, where $p > 0$.
Clearly $L_1(\alpha; m) = \infty$. Furthermore, the function
\begin{equation}
f(x) = x^p
\label{SE13}
\end{equation}
satisfies the conditions of (\ref{C1}), thus Theorems 2.3--2.5 can be applied. Formula \eqref{delta} becomes
\begin{equation}
\delta = \frac{1}{\ln \left(\frac{f(N)}{f^{\prime}(N)}\right)} = \frac{1}{\ln N - \ln p}.
\label{SE12}
\end{equation}
Hence, Theorems 2.3--2.5 yield
\begin{equation}
E\left[\, T_m(N)\,\right] = A_{N} N^p \left[\,\ln N + (m-2)\ln\ln N - \ln p + \gamma -\ln\left(m-1\right)! + o(1) \, \right]
\label{SE14}
\end{equation}
\begin{equation}
E\left[\,T_m(N) \left(T_{m}(N)+1\right)\,\right] = A^2_N N^{2p} (\ln N)^2 \left[1 + O\left(\frac{\ln\ln N}{\ln N}\right)\right],
\label{SE15}
\end{equation}
and
\begin{equation}
V\left[\,T_m(N)\,\right] \sim \frac{\pi^2}{6} \, A^2_N \, N^{2p}
\label{SE16}
\end{equation}
where
\begin{equation}
A_N = \sum_{j=1}^N \frac{1}{j^p}.
\label{SE17}
\end{equation}
Thus (see, e.g., \cite{A}),
\begin{equation}
A_N = \frac{N^{1-p}}{1-p} + \zeta(p) + O\left(\frac{1}{N^p}\right),
\qquad \text{if }\; 0 < p < 1,
\label{SE18}
\end{equation}
\begin{equation}
A_N = \ln N + \gamma + O\left(\frac{1}{N}\right),
\qquad \text{if }\; p = 1,
\label{SE19}
\end{equation}
and
\begin{equation}
A_N = \zeta(p) - \frac{1}{(p-1) N^{p-1}} + O\left(\frac{1}{N^p}\right),
\qquad \text{if }\; p > 1,
\label{SE20}
\end{equation}
where $\zeta(\cdot)$ is Riemann's Zeta function (recall that $\zeta(p) < 0$ if $0 < p < 1$). For instance, if $0 < p < 1$, then, in view of
\eqref{SE18}, formulas \eqref{SE14} and \eqref{SE16} yield
\begin{equation}
E\left[\, T_m(N)\,\right] = \frac{N}{1-p} \, \left[\,\ln N + (m-2)\ln\ln N - \ln p + \gamma -\ln\left(m-1\right)! \, \right]
+ O\left(N^p \ln N\right)
\label{SE21}
\end{equation}
and
\begin{equation}
V\left[\,T_m(N)\,\right] \sim \frac{\pi^2}{6} \, \frac{N^2}{(1-p)^2}
\label{SE22}
\end{equation}
respectively. Formula \eqref{SE21} should be compared with \eqref{1}--\eqref{2}; likewise formula \eqref{SE22} should be related to
formula \eqref{SE22c} of Conjecture 1.

\smallskip

\textbf{Example 3} (the exponential law)\textbf{.} As in Remark 2, let $p > 0$ and consider the sequences $\alpha = \{e^{pj} \}_{j=1}^{\infty}$
and $\beta = \{e^{-pj} \}_{j=1}^{\infty}$. We have already observed that $\alpha$ and $\beta$ produce the same coupon probabilities.
We have $L_1(\beta; m) = \infty$. Furthermore $f(x) = e^{px}$ does not satisfy condition (ii) of (\ref{C1}), thus Theorems 2.3--2.5 cannot be
applied. Let us consider, instead, the sequence $\alpha$, where we have $L_1(\alpha; m) < \infty$. Here
\begin{equation}
A_N = \sum_{j=1}^N e^{pj} = \frac{e^{p(N+1)} - e^p}{e^p - 1} = \frac{e^{p(N+1)}}{e^p - 1} + O(1),
\qquad
N \to \infty.
\label{SE4}
\end{equation}
Also, formulas \eqref{ET2} and \eqref{ET3} of Remark 3 give
\begin{equation}
\delta _{N} = O\left(e^{-pN}\right),
\qquad \text{and}\qquad
\Delta _{N} = O\left(e^{-2pN}\right).
\label{SE6}
\end{equation}
Therefore, Theorem 2.2 yields (as $N \to \infty$)
\begin {align}
E\left[\,T_{m}(N)\,\right] &= \frac{e^{p(N+1)}}{e^p - 1} L_{1}(\alpha;m ) + O(1),
\label{SE9}
\\
E\left[\,T_{m}(N)(T_{m}(N)+1)\,\right] &= \frac{e^{2p(N+1)}}{(e^p - 1)^2} L_{2}(\alpha; m) + O\left(e^{pN}\right),
\label{SE10}
\end{align}
and
\begin{equation}
V\left[ T_{m}(N)\right] = \frac{e^{2p(N+1)}}{(e^p - 1)^2} \left(L_2(\alpha; m) - L_1(\alpha; m)^2\right)
+ O\left(e^{pN}\right).
\label{SE11}
\end{equation}
It follows that, regarding the sequence $\beta$, the associated asymptotics are also given by \eqref{SE9}, \eqref{SE10}, and \eqref{SE11}.
In this way we get cheaply a counterexample for Theorems 2.3--2.5, in case where $f(\cdot)$ does not satisfy all conditions of (\ref{C1}).

\section{Limit Distributions}

\subsection{Case I: $L_1(\protect\alpha; m) < \infty$}

\smallskip

\begin{theorem}
Let $\alpha =\{a_j\}_{j=1}^{\infty}$ be a sequence such that $L_1(\protect\alpha;m ) < \infty$
(recall (\ref{17}) and Theorem 2.1) and, as in Section 2,
\begin{equation*}
p_j = \frac{a_j}{A_N},
\qquad \text{where}\quad
A_N = \sum_{j=1}^N a_j.
\end{equation*}
Then, for all $s \in [0, \infty)$ we have
\begin{equation}
P\left\{\frac{T_m(N)}{A_N} \leq s \right\}\rightarrow F(s) := \prod_{j=1}^{\infty}\left[1 - S_{m}(a_j s) e^{-a_j s} \right],
\qquad
N \rightarrow \infty,
\label{PR0}
\end{equation}
where $S_m(\,\cdot \,)$ is given by (\ref{7}).
\end{theorem}
\textit{Proof.} Setting $z = e^{\lambda}$ with $\Re(\lambda) > 0$, formula (\ref{P8}) can be written as
\begin{equation}
E\left[ e^{-\lambda T_m(N)} \right]
= 1 - (e^{\lambda} - 1) \int_0^{\infty}
\left\{1-\prod_{j=1}^{N}\left[1 - S_{m}(p_j t) e^{-p_j t} \right]\right\} \exp\left(-(e^{\lambda} - 1) t \right) dt,
\label{PR1}
\end{equation}
where $\Re(\lambda) > 0$.

Substituting $t = A_N s$ in the integral of (\ref{PR1}) we obtain
\begin{equation*}
E\left[ e^{-\lambda T_m(N)} \right]
= 1 - (e^{\lambda} - 1) A_N
\int_0^{\infty} \left\{1-\prod_{j=1}^{N}\left[1 - S_{m}(a_j s) e^{-a_j s} \right]\right\} \exp\left(-(e^{\lambda} - 1) A_N s \right) ds,
\end{equation*}
or
\begin{equation*}
E\left[ e^{-\lambda T_m(N) / A_N} \right] =
\qquad\qquad\qquad\qquad\qquad\qquad\qquad\qquad\qquad\qquad\qquad\qquad\qquad\qquad
\end{equation*}
\begin{equation}
1 - (e^{\lambda / A_N} - 1) A_N
\int_0^{\infty} \left\{1-\prod_{j=1}^{N}\left[1 - S_{m}(a_j s) e^{-a_j s} \right]\right\} \exp\left(-(e^{\lambda / A_N} - 1) A_N s \right) ds.
\label{PR2}
\end{equation}
Finally, in view of Subsection 2.2 and (\ref{PR0}) (for the definition of $F$) dominated convergence gives
\begin{equation*}
\lim_N E\left[ e^{-\lambda T_m(N) / A_N} \right]
= 1 - \lambda
\int_0^{\infty} \left[1 - F(s) \right] e^{-\lambda s} ds = \int_0^{\infty} e^{-\lambda s} dF(s),
\end{equation*}
for all complex $\lambda$ such that $\Re(\lambda) > 0$.
\hfill $\blacksquare$

\smallskip

Notice that the limit distribution depends on the sequence $\alpha$.

\subsection{Case II: $L_1(\protect\alpha; m) = \infty$}

P. Neal \cite{N} has established a general theorem regarding the limit distribution of $T_m(N)$ (appropriately normalized) as $N \to \infty$,
where $\pi_N = \{p_{N1}, p_{N2},...,p_{NN} \}$, $N = 1, 2,...$, is a sequence of (sub)probability measures, not necessarily of the form (\ref{8}).

\smallskip

\textbf{Theorem N.} Suppose that there exist sequences $\{b_N\}$ and $\{k_N\}$ such that $k_N / b_N \rightarrow 0$
as $N \rightarrow \infty$ and that, for $y \in \mathbb{R}$,
\begin{equation}
\Lambda_N(y\,;m) := \frac{b_N^{m-1}}{\left(m-1\right)!} \sum_{j=1}^N p_{Nj}^{m-1}\exp\bigg(-p_{Nj} \left(b_N + y k_N\right) \bigg) \rightarrow g(y),
\quad
N \rightarrow \infty,
\label{N1}
\end{equation}
for a nonincreasing function $g(\cdot)$ with $g(y) \rightarrow \infty$ as $y \rightarrow -\infty$ and
$g(y) \rightarrow 0$ as $y \rightarrow \infty$. Then
\begin{equation}
\frac{T_{m}(N) - b_N}{k_N} \overset{D}{\longrightarrow} Y,
\qquad
N \rightarrow \infty,
\label{N2}
\end{equation}
where $Y$ has distribution function
\begin{equation}
F(y) = P\{ Y \leq y \} = e^{-g(y)},
\qquad
y \in \mathbb{R}.
\label{N222a}
\end{equation}.

Theorem N \textit{does not indicate at all} how to choose the sequences $\{b_N\}$ and $\{k_N\}$. Here our asymptotic formulas can help.

The conclusion (\ref{N2}) of Theorem N suggests that as $N \rightarrow \infty$
\begin{equation*}
b_N \sim E[T_{m}(N)] \qquad \text{and} \qquad k_N \sim c\sqrt{V[T_{m}(N)]},\quad \text{for some}\ c \neq 0.
\end{equation*}
Recall that for the Case II the coupon probabilities $p_{Nj}$, $1 \leq j \leq N$, $N = 1,2,...$, are taken as
\begin{equation}
p_{Nj} = \frac{a_j}{A_N} \qquad \text{with} \qquad A_N = \sum_{j=1}^N a_j,
\quad
a_j = \frac{1}{f(j)},
\label{N2a}
\end{equation}
where $f(x)$ satisfies (\ref{C1}). Then Theorems 2.3 and 2.5 propose the choices
\begin{equation}
b_N = A_N f(N) \big[\rho(N) + \left(m-2\right)\ln \rho(N)\big] \qquad  \text{and} \qquad k_N = A_N f(N),
\label{N6}
\end{equation}
where
\begin{equation}
\rho(N) := 1 / \delta = \ln \left(f(N) / f'(N)\right)
\label{N6aa}
\end{equation}
(notice that, as $N \rightarrow \infty$ we have that $\rho(N) \rightarrow \infty $, and hence $k_N / b_N \rightarrow 0$ as required).
In this case, $\Lambda_N(y\,;m)$ of (\ref{N1}) becomes
\begin{equation}
\Lambda_N(y\,;m) = \frac{f(N)^{m-1}}{\left(m-1\right)!} \left[\rho(N) + \left(m-2\right)\ln \rho(N)\right]^{m-1} \;\tilde{\Lambda}_N(y\,;m),
\label{C1bb}
\end{equation}
where
\begin{align*}
\tilde{\Lambda}_N(y\,;m):=\sum_{j=1}^N \frac{1}{f\left(j\right)^{m-1}}\exp\left(-\frac{f(N)}{f(j)} \left[\rho(N) + \left(m-2\right)\ln \rho(N)+y\right] \right).
\end{align*}
Since $f$ is increasing and satisfies (\ref{C1}) we have for sufficiently large $N$
\begin{align}
\tilde{\Lambda}_N(y\,;m)&=\int_{1}^N \frac{1}{f\left(x\right)^{m-1}}\exp\left(-\frac{f(N)}{f(x)} \left[\rho(N) + \left(m-2\right)\ln \rho(N)+y\right] \right)\nonumber\\
&+O\left(\frac{1}{f\left(N\right)^{m-1}}\exp\left(-\left[\rho(N) + \left(m-2\right)\ln \rho(N)+y\right] \right)\right).
\label{error}
\end{align}
Let us consider the integral
\begin{equation*}
\tilde{I}_N(y\,;m): =\int_{1}^N \frac{1}{f\left(x\right)^{m-1}}\exp\left(-\frac{f(N)}{f(x)} \left[\rho(N) + \left(m-2\right)\ln \rho(N)+y\right] \right).
\end{equation*}
Integration by parts gives
\begin{align}
\tilde{I}_N(y\,;m) =& \left[ \frac{1}{M} \cdot \frac{f(x)^{3-m}}{f'(x)} \exp\left(-\frac{M}{f(x)}\right) \right]_{x=1}^N \nonumber \\
&+\int_1^N \frac{f(x)^{2-m}}{M}\left[m-3+\frac{f^{\prime\prime}(x)/f^{\prime}(x)}{f^{\prime}(x)/f(x)} \right] \exp\left(-\frac{M}{f(x)}\right)dx, \label{N5}
\end{align}
where for typographical convenience we have set
\begin{equation*}
M := f(N) \left[\rho(N) + \left(m-2\right)\ln \rho(N)+y\right].
\end{equation*}
By (\ref{C1}) it follows that the integral in the right-hand side of (\ref{N5}) is $o(\tilde{I}_N(y\,;m))$ as $N\rightarrow \infty$. The quantity
$\frac{1}{M} \cdot \frac{f(1)^{3-m}}{f'(1)} \exp\left(-\frac{M}{f(1)}\right)$ is, also, $o(\tilde{I}_N(y\,;m))$. Hence, as $N\rightarrow \infty$
\begin{equation*}
\tilde{\Lambda}_N(y\,;m) \sim \tilde{I}_N(y\,;m) \sim \frac{f(N)^{2-m}}{f'(N)} \cdot \frac{\exp\left(-\rho(N) -\left(m-2\right) \ln \rho(N) - y \right)}{\rho(N) +\left(m-2\right) \ln \rho(N) + y}.
\end{equation*}
In view of (\ref{N6aa}) and the fact that $\rho(N) \to \infty$ as $N \to \infty$, the above formula becomes
\begin{equation*}
\tilde{\Lambda}_N(y\,;m) \sim f(N)^{1-m} \rho(N)^{1-m} e^{-y}.
\end{equation*}
Using the above asymptotics in (\ref{C1bb}) yields
\begin{equation*}
\Lambda_N(y\,;m) \rightarrow \frac{e^{-y}}{(m-1)!},
\quad
N \rightarrow \infty.
\end{equation*}
Therefore, by invoking Theorem $N$ we obtain the following limit theorem.

\smallskip

\begin{theorem}
Suppose the coupon probabilities come from a sequence $\alpha$ as in (\ref{N2a}), where $f(x)$ satisfies (\ref{C1}). Then, for all $y \in \mathbb{R}$ we have
\begin{equation}
P\left\{\frac{T_m(N) - b_N}{k_N} \leq y \right\}\rightarrow \exp\left(-\frac{e^{-y}}{\left(m-1\right)!}\right),
\qquad
N \rightarrow \infty,
\label{N7}
\end{equation}
where $b_N$ and $k_N$ are given by (\ref{N6})--(\ref{N6aa}).
\end{theorem}
Notice that the limiting distribution in (\ref{N7}) is \textit{independent} of the choice of $f(x)$.

\smallskip

\textbf{Remark 7.} The fact that for the sequences $b_N$ and $k_N$ of (\ref{N6})--(\ref{N6aa}) the limit $g(y)$ in (\ref{N1}) exists
and has the right behavior is an independent confirmation that the statements of the Theorems 2.3, 2.4, and 2.5 are correct.

\smallskip

\textbf{Example 2} (continued)\textbf{.} For $p > 0$ let us take $f(x) = x^p$, so that $a_j = 1 / j^p$. Then, \eqref{N6aa} becomes
\begin{equation}
\rho(N) = \ln N - \ln p
\label{N6aaee}
\end{equation}
and hence $b_N$ and $k_N$ of \eqref{N6} can be taken as
\begin{equation}
b_N = A_N N^p \big[\ln N + (m-2)\ln\ln N - \ln p \big]
\qquad  \text{and} \qquad
k_N = A_N N^p,
\label{N6ee}
\end{equation}
where $A_N$ is given by \eqref{SE18}--\eqref{SE20} (it is enough to use the leading asymptotic term of $A_N$). If, in particular,
$0 < p < 1$, then formula \eqref{N7} holds with
\begin{equation}
b_N = \frac{N}{1-p} \big[\ln N + (m-2)\ln\ln N - \ln p \big]
\qquad  \text{and} \qquad
k_N = \frac{N}{1-p}.
\label{N6eee}
\end{equation}
This example should be compared with the limiting behavior \eqref{3} of the case of equal coupon probabilities.

\subsection{Slowly decaying sequences}
Suppose that our sequence $\alpha = \{a_j\}$ decays to $0$ slower that $1 / j^p$ for every $p > 0$ (of course, $L_1(\protect\alpha; m) = \infty$).
Then, the corresponding function $f(x)$ (for which $a_j = 1 / f(j)$) may not satisfy all conditions of \eqref{C1} and, hence, the method presented in
Subsection 2.4 for determining the asymptotics of the expectation and the variance of $T_{m}(N)$ may not work. Nevertheless,
Example 2 together with formulas \eqref{3}--\eqref{333} for the case of equal coupon probabilities (i.e. when $a_j = \;$constant) suggest that the
limit distribution of $T_{m}(N)$, appropriately normalized, should be Gumbel and, furthermore that the sequences $b_N$ and $k_N$ of Theorem N should
be taken as
\begin{equation}
b_N = N \ln N + (m-1) N \ln\ln N
\qquad  \text{and} \qquad
k_N = N
\label{SD0}
\end{equation}
(also, that, as $N \to \infty$, $V[T_{m}(N)] \sim (\pi^2 / 6) N^2$, while $E[T_{m}(N)] = N \ln N + (m-1) N \ln\ln N + c_m(\alpha) N + o(N)$, where $c_m(\alpha)$ is a constant depending on $m$ and $\alpha$ such that $c_m(\alpha) \geq \gamma - \ln(m-1)!$).

Let us illustrate the above comment with the function $f(x) = (\ln x)^p$, $p > 0$, which does not satisfy conditions (iii) and (iv) of \eqref{C1}:

\smallskip

Suppose the coupon probabilities come from the sequence $\alpha = \{a_j = (\ln j)^{-p}\}_{j=2}^{\infty}$ for some $p > 0$.
Then, for all $y \in \mathbb{R}$ we have
\begin{equation}
P\left\{\frac{T_m(N) - N \ln N - (m-1) N \ln\ln N}{N} \leq y \right\}\rightarrow \exp\left(-\frac{e^{-(y-p)}}{(p+1) (m-1)!}\right)
\label{SD1}
\end{equation}
as $N \rightarrow \infty$ (the proof of formula \eqref{SD1} can be found in \cite{DPL}). Needless to say that formula \eqref{SD1} is equivalent to
\begin{equation}
P\left\{\frac{T_m(N) - N \ln N - (m-1) N \ln\ln N - \left[\gamma + p - \ln(p+1) - \ln(m-1)! \right] N }{N} \leq y \right\}\rightarrow e^{-e^{-y}}
\label{SD1a}
\end{equation}
as $N \rightarrow \infty$.

Finally, let us observe that for $a_j = (\ln j)^{-p}$, $p > 0$, the above example suggests the asymptotic formulas (as $N \rightarrow \infty$)
\begin{equation}
E\left[T_m(N)\right] = N \ln N + (m-1) N \ln\ln N + \left[\gamma + p - \ln(p+1) - \ln(m-1)! \right] N + o(N),
\label{SD0a}
\end{equation}
and
\begin{equation*}
V\left[T_m(N)\right] \sim \frac{\pi^2}{6} \, N^2.
%\label{SD0b}
\end{equation*}
Notice that the expected value in \eqref{SD0a} is slightly bigger than the corresponding expected value for the case of equal coupon probabilities
(recall \eqref{1}--\eqref{2}), due to the term $p - \ln(p+1)$ which is strictly positive for all $p > 0$.

\section{Appendix}
Here we give the proofs of Theorem 2.1 and some technical lemmas which appeared in Section 2.\\\\
{\it Proof of Theorem 2.1}. Before proving the theorem we recall the following inequality which can be proved easily by induction and limit:\\
Let $\left\{ b_{j}\right\} _{j=1}^{\infty }$ be a sequence
of real numbers such that $0\leq b_{j}\leq 1,$ for all $j$.
If $\sum_{j=1}^{\infty }b_{j}<\infty $, then
\begin{equation}
1-\prod_{j=1}^{\infty }(1-b_{j})\leq \sum_{j=1}^{\infty }b_{j}.
\label{LL1b}
\end{equation}
Let us prove the equivalence of $(i)$ and $(iii)$. The equivalence between $(ii)$ and $(iii)$ is similar.\\
Assume that there is a $\xi \in (0,1)$ such that (\ref{heat})
is true. Then, by (\ref{17}) and (\ref{LL1b}) we have
\begin{equation*}
L_1 (\alpha;m)\leq \int_{0}^{\xi}\sum_{j=1}^{\infty} \big(x^{a_{j}}S_{m}\left(-a_{j}\ln x\right)\big)\frac{dx}{x}
+\int_{\xi }^{1}\frac{dx}{x}.
\end{equation*}
Using Tonelli's Theorem we have, in view of (\ref{7})
\begin{equation*}
L_1 (\alpha;m)\leq\sum_{j=1}^{\infty}\sum_{k=0}^{m-1}\left(\frac{\left(-1\right)^{k}a^{k}_{j}}{k!}\int_{0}^{\xi}x^{a_{j}-1}\left(\ln x\right)^{k}\,dx\right)+ \ln \xi.
\end{equation*}
For the integral above we have by repeated integration by parts
\begin {equation*}
\int_{0}^{\xi}x^{a_{j}-1}\left(\ln x\right)^{k}\,dx=\frac{1}{a_{j}}\;\xi^{a_{j}}\sum_{i=0}^{k}(-1)^{i}\left(k\right)_{i}\frac{1}{a^{i}_{j}}
\left(\ln \xi\right)^{k-i},
\end{equation*}
where $\left(k\right)_{i}=k!/\left(k-i\right)!$ is the falling Pochhammer symbol. Hence,
\begin{equation*}
L_1 (\alpha;m)\leq\sum_{j=1}^{\infty}\left[\sum_{k=0}^{m-1}\frac{\left(-1\right)^{k}a^{k}_{j}}{k!}\left(\frac{1}{a_{j}}\;\xi^{a_{j}}
\sum_{i=0}^{k}(-1)^{i}\left(k\right)_{i}\frac{1}{a^{i}_{j}}\left(\ln \xi\right)^{k-i}\right)\right]+ \ln \xi.
\end{equation*}
Now, (\ref{heat}) implies that $\xi ^{a_{j}}\rightarrow 0$,
hence $a_{j}\rightarrow \infty $. Therefore, $\min_{j}\left\{ a_{j}\right\} = a_{j_{0}} > 0$.
Thus,
\begin{equation*}
L_1 (\alpha;m)\leq\,\left(\sum_{j=1}^{\infty}\xi^{a_{j}}a^{m-1}_{j}\right)\Bigg[\sum_{k=0}^{m-1}\frac{\left(-1\right)^{k}}{k!}\left(\ln \xi\right)^{k}
\left(\sum_{i=0}^{k}(-1)^{i}\left(k\right)_{i}\frac{\left(\ln \xi\right)^{-i}}{a^{i}_{j_{0}}}\right)\Bigg]
+ \ln {\xi}.
\end{equation*}
Since $\xi \in (0,1)$, (\ref{heat}) implies
\begin{equation*}
\sum_{j=1}^{\infty}\xi^{a_{j}}a^{m-1}_{j}<\infty.
\end{equation*}
It follows that $L_1 (\alpha;m)<\infty.$ Conversely, if
\begin{equation*}
\sum_{j=1}^{\infty} \xi^{a_j}=\infty  \quad \text {for \ all }\  \xi \in (0,1)
\end{equation*}
then for any fixed positive integer $m$  we have
\begin{equation*}
\sum_{j=1}^{\infty} \xi^{a_j}a^{m-1}_{j}=\infty , \quad \text {for \ all }\  \xi \in (0,1)
\end{equation*}
and by a standard property of infinite products (see, e.g., \cite{Ru}) it follows that
\begin{equation*}
\prod_{j=1}^{\infty}\bigg(1-x^{a_{j}}\,S_{m}\left(-a_{j}\ln x\right)\bigg) = 0,  \quad \text {for \ all }\ x  \in (0,1).
\end{equation*}
Hence (\ref{17}) yields $L_{1}(\alpha;m )= \int_{0}^{1} (dx / x) = \infty $.
$\hfill \blacksquare$ \\\\
\textit{Proof of Lemma 2.2 - PART I} (the integral $I_1$).\\
Regarding the integral of (\ref{I1}), given any $\varepsilon \in (0,1)$ we have
\begin{align*}
I_{1}(N) :&= \int_0^{1 - \varepsilon} \exp \left[
\sum_{j=1}^{N}\ln \left[ 1-e^{-\frac{F(N)}{f(j)}s}S_{m}\left(\frac{F(N)}{f(j)}s\right)\right] \right] ds\nonumber\\
&< \exp \left[-\sum_{j=1}^N e^{-\frac{F(N)}{f(j)} \left(1 - \varepsilon \right)} S_{m}\left(\frac{F(N)}{f(j)}\left(1-\varepsilon\right)\right)\right]\nonumber\\
&< \exp \left[-\sum_{k=0}^{m-1}\left[\frac{\left(1-\varepsilon\right)^{k}F(N)^{k}}{\left(m-1\right)!}
\left( \sum_{j=1}^{N} f(j)^{-k}e^{-\left(1 - \varepsilon \right)\frac{F(N)}{f(j)}}\right)\right]\right],
\end{align*}
since $\ln(1-x) < -x$, for $0 < x < 1$. Let us now consider the function
\begin{equation*}
g(x):= f(x)^{-k} \exp\left(-\lambda F(N)/f(x)\right),\,\,\,x\in [1,N],\,\,\,k=0,1,\cdots,m-1,\,\,\,\,\,\, \lambda\in(0,1).
\end{equation*}
It is easy to check that conditions (\ref{C1}) imply that for sufficiently large $N$ $g(\cdot)$ is strictly increasing in $[1,N].$ Hence,
\begin{align*}
\int_{1}^{N}g(x)\,dx&\leq \sum_{j=1}^{N}g(j)\leq \int_{1}^{N}g(x)\,dx +g\left(N\right).
\end{align*}
Moreover, by Lemma 2.1 it is easy to see that $g(N)=o\left(\int_{1}^{N}g(x)\,dx\right)$ as $N\rightarrow \infty$. Thus,
\begin{equation*}
\sum_{j=1}^{N}g(j)\sim \int_{1}^{N}g(x)\,dx \quad \text{as}\quad N\rightarrow \infty.
\end{equation*}
Applying Lemma 2.1 for $\kappa =-k$ one arrives at
\begin{equation*}
I_{1}(N)
< \exp\left[-\sum_{k=0}^{m-1}\frac{\left(1-\varepsilon\right)^{k}F(N)^{k}}{\left(m-1\right)!}\left[\frac{f(N)^{2-k}}{(1 - \varepsilon) F(N) f'(N)} e^{-\frac{F(N)}{f(N)}(1 - \varepsilon)}\left(1 + M_1\;\frac{f(N)}{F(N)} \right) \right]\right],
\end{equation*}
where $M_1$ is a positive constant.
Using (\ref{F}) and (\ref{delta}) i.e. the definitions of $F(\cdot)$ and $\delta$, we have
\begin{align*}
I_{1}(N)&
< \exp\left[-\sum_{k=0}^{m-1}\frac{\left(1-\varepsilon\right)^{k-1}}{\left(m-1\right)!}\,\frac{e^{\varepsilon/\delta}}{\delta^{k-1}}\,\left(1 + M_1\,\delta \right) \right]\nonumber\\
&=\exp\left[-\left(1-\varepsilon\right)\frac{\delta^{m}-\left(1-\varepsilon\right)^{m}}{\delta^{m}\left(\delta-\left(1-\varepsilon\right)\right)}\,e^{\varepsilon/\delta}\,\left(1 + M_1\,\delta \right) \right].
\end{align*}
Since $\delta \rightarrow 0^{+}$ and $\varepsilon \in (0,1)$ we have
\begin{equation*}
I_{1}(N) << \delta^{4-m}e^{-\varepsilon/ \delta},
\end{equation*}
for sufficiently large $N$, $m=1,2,3,\cdots.$ $\hfill \blacksquare$ \\\\
\textit{Proof of Lemma 2.2 - PART II} (the integral $I_2$).\\
Our first task is to compute a few terms of the asymptotic expansion of the integral of (\ref{I2}). For convenience we set
\begin{equation}
B_{m}(N;s) :=\sum_{j=1}^N \ln \left[ 1-e^{-\frac{F(N)}{f(j)}s}S_{m}\left(\frac{F(N)}{f(j)}s\right)\right].
\label{AsN}
\end{equation}
Since
\begin{equation*}
\frac{F(N)}{f(j)}\rightarrow \infty \qquad \text{as }\ N\rightarrow
\infty,
\end{equation*}
and $\ln (1-x) = -x + O(x^{2})$ as $x\rightarrow 0$, we have (as long as $s\geq s_{0} > 0$)
\begin{equation}
B_{m}(N;s)=\sum_{j=1}^{N} \left[-e^{-\frac{F(N)}{f(j)}s}S_{m}\left(\frac{F(N)}{f(j)}s\right)+O\left(e^{-\frac{2F(N)}{f(j)}s}\left[S_{m}\left(\frac{F(N)}{f(j)}s\right)\right]^{2}\right)\right].\label{O1}
\end{equation}
Using (\ref{7}), (\ref{O1}) yields
\begin{align}
B_{m}(N;s)=&-\sum_{k=0}^{m-1}\frac{F(N)^{k}\,s^{k}}{k!}
\left(\sum_{j=1}^{N}f(j)^{-k}\,e^{-\frac{F(N)}{f(j)}s}\right) \nonumber \\
&+\sum_{j=1}^{N}O\left(e^{-\frac{2F(N)}{f(j)}s}\left[S_{m}\left(\frac{F(N)}{f(j)}s\right)\right]^{2}\right).\label{40}
\end{align}
Since $f(\cdot)$ is increasing and under conditions (\ref{C1}), it follows from the comparison of sums and integrals that for sufficiently large  $N$
\begin{equation}
\sum_{j=1}^{N}f(j)^{-k}\,e^{-\frac{F(N)}{f(j)}s}=\int_{1}^{N}f(x)^{-k}\,e^{-\frac{F(N)}{f(x)}s}dx+O\left(f(N)^{-k}\,e^{-\frac{F(N)}{f(N)}s}\right).\label{41}
\end{equation}
In view of (\ref{41}) and Lemma 2.1 (for $\kappa = -k$), (\ref{41}) yields (as long as $s\geq s_{0}>0$),
\begin{align}
B_{m}(N;s)=-\sum_{k=0}^{m-1}\frac{F(N)^{k}\,s^{k}}{k!}&\left[\frac{f(N)^{2-k}}{sF(N)f^{\prime}(N)}e^{-\frac{F(N)}{f(N)}s}\right.\nonumber\\
&\left.+\omega(N)\;\frac{f(N)^{3-k}}{s^{2}F(N)^{2}f^{\prime}(N)}e^{-\frac{F(N)}{f(N)}s}\left[1+O\left(\frac{f(N)}{F(N)}\right)\right]\right].\label{42}
\end{align}
For typographical convenience we set
\begin{equation}
A := \frac{f(N)}{f'(N)}.
\label{A}
\end{equation}
(notice that $A \rightarrow \infty$ as $N \rightarrow \infty$). Using (\ref{F}) and (\ref{A}), (\ref{42}) yields
\begin{align}
B_{m}(N;s)=&-\frac{1}{\left(m-1\right)!}A^{1-s}s^{m-2}\left(\ln A\right)^{m-2}\nonumber \\
&-\left(\omega(N)+\frac{1}{\left(m-2\right)!}\right)A^{1-s}s^{m-3}\left(\ln A\right)^{m-3}\left[1+O\left(\frac{1}{s\ln A}\right)\right].
\label{OO2}
\end{align}
Hence the quantity (see, (\ref{I2}))
\begin{equation*}
I_{2}(N) := \int_{1-\varepsilon}^{1+\varepsilon} e^{B_{m}(N;s)}ds
\end{equation*}
via the substitutions $s = 1 - t$ and $u = A^t \left(\ln A\right)^{m-2}$ (and in view of (\ref{7})), yields
\begin{align*}
I_{2}(N)=&\delta\int_{\delta^{2-m}\exp\left(-\varepsilon /\delta\right)}^{\delta^{2-m}\exp\left(\varepsilon /\delta\right)}
\exp\left\{-\frac{1}{\left(m-1\right)!}\,u\left[\,1-\delta \ln u-\left(m-2\right)\delta\ln \delta\,\right]^{m-2}\right.\nonumber\\
&\left.-\left(\omega(N)+\frac{1}{\left(m-2\right)!}\right)\,u\,\delta\,\left[\,1-\delta \ln u-\left(m-2\right)\delta\ln \delta\,\right]^{m-3}
\left(1+O\left(\delta\right)\right)\right\}\frac{du}{u},\nonumber
\end{align*}
where (see (\ref{delta}))
\begin{equation*}
 \delta := \frac{1}{\ln A}=\frac{1}{\ln \left(\frac{f(N)}{f^{\prime}(N)}\right)}=\frac{f(N)}{F(N)}
\end{equation*}
(hence, $A\rightarrow \infty$ implies $\delta \rightarrow 0^+$). We have
\begin{equation}
I_{2}=\delta \left(\int_{\delta^{2-m}\exp\left(-\varepsilon /\delta\right)}^{1/\sqrt{\delta}}
+ \int_{1/\sqrt{\delta}}^{\delta^{2-m}\exp\left(\varepsilon /\delta\right)}\right). \label{eint5}
\end{equation}
First we get an upper bound for the second integral of (\ref{eint5}) as follows:
\begin{align}
\int_{1/\sqrt{\delta}}^{\delta^{2-m}\exp\left(\varepsilon /\delta\right)}
\exp&\left\{-\frac{1}{\left(m-1\right)!}\,u \left(1-\delta \ln\left(u\,\delta^{m-2}\right)\right)^{m-2}\right.\nonumber\\
&\left.\left[ 1 + \left(\omega(N)+\frac{1}{\left(m-2\right)!}\right)\frac{\delta\left( 1 + O\left(\delta\right)\right)}{1-\delta \ln\left(u\,\delta^{m-2}\right)}\;\right]\right\}\frac{du}{u},\nonumber\\
=O&\left(\sqrt{\delta}\; e^{-1/\left(m-1\right)!\sqrt{\delta}}\right).\label{eint6}
\end{align}
Let us denote $K_1(\delta)$ the first integral of (\ref{eint5}). We use the binomial theorem to expand the quantities $\left[1-\delta \ln u-\left(m-2\right)\delta\ln \delta\right]^{m-2}$ and $\left[1-\delta \ln u-\left(m-2\right)\delta\ln \delta\right]^{m-3}$. Next, we expand the exponentials and get
\begin{align*}
K_{1}(\delta)=\int_{\delta^{2-m}\exp\left(-\varepsilon /\delta\right)}^{1/\sqrt{\delta}}\frac{e^{-u/\left(m-1\right)!}}{u}\;
&\left\{1+\frac{m-2}{\left(m-1\right)!}u\,\delta \ln\left(u\,\delta^{m-2}\right)+u^{2}O\left[\,\delta \ln\left(u\,\delta^{m-2}\right)\right]^{2}\right\}\nonumber\\
\times&\left\{1-\left(\omega(N)+\frac{1}{\left(m-2\right)!}\right)\,u\,\delta\,\left(1+O\left(\delta\right)\right)\right.\nonumber\\
&\left.+u^{2}O\left[\,\delta \ln\left(u\,\delta^{m-2}\right)\right]^{2}\right\}{du}
\end{align*}
(since $e^x =1 + x + O(x^{2})$ as $x \rightarrow 0$). Hence,
\begin{align*}
K_{1}(\delta)=\int_{\delta^{2-m}\exp\left(-\varepsilon /\delta\right)}^{1/\sqrt{\delta}}\frac{e^{-u/\left(m-1\right)!}}{u}
&\left[1+\frac{m-2}{\left(m-1\right)!}\,u\delta\,\ln\left(u\,\delta^{m-2}\right)\right.\nonumber\\
&\left.\,\,\,\,\,\,\,-\left(\omega(N)+\frac{1}{\left(m-2\right)!}\right)\,u\,\delta\,\left(1+O\left(\delta\right)\right)\right.\nonumber\\
&\left.\,\,\,\,\,\,\,+u^{2}O\left[\delta \ln\left(u\,\delta^{m-2}\right)\right]^{2}\right]du.
\end{align*}
We split the integral above as
\begin{equation}
K_{1}(\delta)=\int_{\delta^{2-m}\exp\left(-\varepsilon /\delta\right)}^{\infty}-\int_{1/\sqrt{\delta}}^{\infty}.\label{20}
\end{equation}
However (and this is an easy exercise)
\begin{align}
&\int_{1/\sqrt{\delta}}^{\infty}\frac{e^{-u/\left(m-1\right)!}}{u}\left[1+\frac{m-2}{\left(m-1\right)!}\,u\delta\,\ln\left(u\,\delta^{m-2}\right)\right.\nonumber\\
&\left.\,\,\,\,\,\,\,\,\,\,\,\,\,\,\,\,\,\,\,\,\,\,\,\,\,\,\,\,\,\,\,\,\,\,\,\,\,\,\,\,\,\,\,\,\,\,\,\,\,\,\,-\left(\omega(N)+\frac{1}{\left(m-2\right)!}\right)\,u\,\delta\,\left(1+O\left(\delta\right)\right)
+u^{2}O\left[\delta \ln\left(u\,\delta^{m-2}\right)\right]^{2}\right]du\nonumber\\
=& O\left(\sqrt{\delta}\,e^{-1/\left(m-1\right)!\sqrt{\delta}}\right) \qquad \text{as $\delta\rightarrow 0^{+}$.}
\label{O}
\end{align}
It follows that in the expression for $K_1(\delta)$ we can replace the upper limit of the integral by $\infty$ and therefore as $\delta\rightarrow 0^{+}$
\begin{align}
I_{2}(N)=\delta\int_{\delta^{2-m}\exp\left(-\varepsilon /\delta\right)}^{\infty}\frac{e^{-u/\left(m-1\right)!}}{u}
&\left[1+\frac{m-2}{\left(m-1\right)!}\,u\delta\,\ln\left(u\,\delta^{m-2}\right)\right.\nonumber\\
&\left.\,\,\,\,\,\,\,-\left(\omega(N)+\frac{1}{\left(m-2\right)!}\right)\,u\,\delta\,\left(1+O\left(\delta\right)\right)\right.\nonumber\\
&\left.\,\,\,\,\,\,\,+u^{2}O\left[\delta \ln\left(u\,\delta^{m-2}\right)\right]^{2}\right]du.\label{eint8}
\end{align}
The following asymptotic expansions easy exercises:
\begin{align}
\int_{x}^{\infty}\frac{e^{-t}}{t}dt
=& -\ln x-\gamma+x+O\left(x^{2}\right)\qquad & \text{as}\quad x\rightarrow 0^{+}, \label{E1}\\
\int_{x}^{\infty}\ln t \;e^{-t}dt
=& -\gamma-x\ln x+x+O\left(x^{2}\,\ln x\right)\qquad & \text{as}\quad x\rightarrow 0^{+}, \label{G11}
\end{align}
where $\gamma =0.5772...$ is the Euler--Mascheroni constant.
Applying (\ref{E1}) and (\ref{G11}) in (\ref{eint8}) we get
\begin{align*}
I_{2}(N)=&\,\varepsilon+\left(m-2\right)\delta \ln \delta+\left[\,\ln\left(m-1\right)!- \gamma\,\right] \delta
+\left(m-2\right)^{2} \delta^{2}\ln \delta\nonumber\\
&+\left[\left(m-2\right)\ln\left(m-1\right)!
-\left(m-2\right)\gamma-\left(m-1\right)-\omega(N)\left(m-1\right)!\right]\delta^{2}\nonumber\\
&+O\left(\delta^{3}\left(\ln \delta\right)^{2}\right).
\end{align*}
Notice that the error term in the above dominates the terms of (\ref{eint6}) and (\ref{O}).$\hfill \blacksquare$ \\\\
\textit{Proof of Lemma 2.2 - PART III} (the integral $I_3$).\\
Our goal is to compute the leading term of $I_{3}(N)$. Here we will follow a different approach.\\
Given $\vartheta \in (0,1)$, there is a $\eta=\eta(\vartheta)$ such that, for $0 < x < \eta$, we have
\begin{equation}
-\left(1+\vartheta\right)x<\ln\left(1-x\right)<-\left(1-\vartheta\right)x
\label{ln}
\end{equation}
and
\begin{equation}
\left(1 - \vartheta\right) x < 1 - e^{-x} < \left(1 + \vartheta\right) x
\label{exp}.
\end{equation}
For $j = 1,\ldots,N$, $s \geq 1$, we use the definition of $F,$ conditions (\ref{C1}), and (\ref{7}) to get
\begin{equation*}
0 < x = e^{-\frac{F(N)}{f(j)}s}\, S_{m}\left(\frac{F(N)}{f(j)}s\right)
=e^{-\frac{F(N)}{f(j)}s}\,\sum_{k=0}^{m-1}\frac{1}{k!}\left(\frac{F(N)}{f(j)}s\right)^{k}
\rightarrow 0 \,\,\text{as}\,\,\ N \rightarrow \infty.
\end{equation*}
Hence, for a given $\vartheta \in (0,1)$, there is $N_0 = N_0(\vartheta)$ such that, for $N \geq N_0$, (\ref{ln}) yields
\begin{align*}
-\left(1 + \vartheta\right) e^{-\frac{F(N)}{f(j)} s}\,S_{m}\left(\frac{F(N)}{f(j)}s\right)
&< \ln\left[1-e^{-\frac{F(N)}{f(j)}s}\,S_{m}\left(\frac{F(N)}{f(j)}s\right)\right]\nonumber\\
&< -\left(1-\vartheta\right) e^{-\frac{F(N)}{f(j)} s}\,S_{m}\left(\frac{F(N)}{f(j)}s\right), \,  j = 1,\ldots,N.
\end{align*}
By summing over $j$ and using (\ref{AsN}) we get
\begin{align*}
-\left(1+\vartheta\right)\sum_{j=1}^{N}e^{-\frac{F(N)}{f(j)}s}\,S_{m}\left(\frac{F(N)}{f(j)}s\right)&<B_{m}(N;s)\nonumber\\
&<-\left(1-\vartheta\right)\sum_{j=1}^N e^{-\frac{F(N)}{f(j)} s}\,S_{m}\left(\frac{F(N)}{f(j)}s\right).
\end{align*}
Using (\ref{7}) we have
\begin{equation*}
\sum_{j=1}^{N}e^{-\frac{F(N)}{f(j)}s}\,S_{m}\left(\frac{F(N)}{f(j)}s\right)=\sum_{k=0}^{m-1}\frac{F(N)^{k}\,s^{k}}{k!}
\left[\sum_{j=1}^{N}f(j)^{-k}\,e^{-\frac{F(N)}{f(j)}s}\right]
\end{equation*}
and from the comparison of sums and integrals (see also  (\ref{41})), we arrive at
\begin{align}
-\left(1+\vartheta\right)&
\left[f(N)^{-k}\,e^{-\frac{F(N)}{f(N)}s} + \int_{1}^{N}f(x)^{-k}\,e^{-\frac{F(N)}{f(x)}s}dx\right]\nonumber\\
&< \sum_{j=1}^{N}f(j)^{-k}\,e^{-\frac{F(N)}{f(j)}s}\nonumber\\
&< -\left(1-\vartheta\right)
\int_{1}^{N}f(x)^{-k}\,e^{-\frac{F(N)}{f(x)}s}dx.
\label{bb1}
\end{align}
Hence,
\begin{align}
-\left(1+\vartheta\right)\sum_{k=0}^{m-1}\frac{F(N)^{k}\,s^{k}}{k!}&
\left[f(N)^{-k}\,e^{-\frac{F(N)}{f(N)}s} + \int_{1}^{N}f(x)^{-k}\,e^{-\frac{F(N)}{f(x)}s}dx\right]\nonumber\\
&<B_{m}(N;s)\nonumber\\
&< -\left(1-\vartheta\right)
\sum_{k=0}^{m-1}\frac{F(N)^{k}\,s^{k}}{k!}\left[
\int_{1}^{N}f(x)^{-k}\,e^{-\frac{F(N)}{f(x)}s}dx\right].
\label{100}
\end{align}
Now, by Lemma 2.1 and from (\ref{SL1}) and (\ref{SL2}) we have  $B_{m}(N;s)\rightarrow 0$ as $N\rightarrow \infty$ uniformly in $s\in [1+\varepsilon,\infty)$, for all positive integers $m$. Thus, for given $\vartheta >0$, there exist $N_1 = N_1(\vartheta)$ such that, for $N \geq N_1$, (\ref{exp}) gives
\begin{equation*}
-\left( 1- \vartheta\right) B_{m}(N;s) < 1 - e^{B_{m}(N;s)} < -\left(1 + \vartheta\right)B_{m}(N;s).
\end{equation*}
Therefore (see (\ref{I2}) and (\ref{AsN})),
\begin{equation*}
-\left(1 - \vartheta\right)\int_{1+\varepsilon}^{\infty} B_{m}(N;s)\,ds<I_{3}(N)<-\left(1+\vartheta\right)\int_{1+\varepsilon}^{\infty}B_{m}(N;s)\,ds.
\end{equation*}
Using the bounds of $B(s; N)$ of (\ref{100}) in the above formula we get that for all $N\geq N_2 = \max \{ N_0, N_1 \}$
\begin{equation*}
\left(1 -\vartheta \right)^2 \int_{1+\varepsilon}^{\infty}\,\int_1^N\,\sum_{k=0}^{m-1}\frac{F(N)^{k}\,s^{k}}{k!}\left[f(x)^{-k}\, e^{-\frac{F(N)}{f(x)}s}dx\;\right]ds
\end{equation*}
\begin{equation*}
-\vartheta\left(1 - \vartheta\right) \sum_{k=0}^{m-1}\frac{F(N)^{k}\,}{f(N)^{k}\,k!}\left[\int_{1+\varepsilon}^{\infty} s^{k}\,e^{-\frac{F(N)}{f(N)}s}ds\right]
\end{equation*}
\begin{equation*}
< I_3(N)
\end{equation*}
\begin{equation*}
< \left(1 + \vartheta\right)^2 \int_{1+\varepsilon}^{\infty}\,\int_1^N
\sum_{k=0}^{m-1}\frac{F(N)^{k}\,s^{k}}{k!}\left[f(x)^{-k}\,e^{-\frac{F(N)}{f(x)}s}dx\;\right]ds
\end{equation*}
\begin{equation}
+ \left(1 + \vartheta\right)^2 \sum_{k=0}^{m-1}\frac{F(N)^{k}\,}{f(N)^{k}\,k!}\left[\int_{1+\varepsilon}^{\infty}s^{k}\,e^{-\frac{F(N)}{f(N)}s}ds\right].
\label{b}
\end{equation}
Using Lemma 2.1, for $\kappa =-k$ we have
\begin{equation*}
\int_{1+\varepsilon}^{\infty}s^{k}\left[\int_1^N
f(x)^{-k}\,e^{-\frac{F(N)}{f(x)}s}dx\;\right]ds=
\frac{f(N)^{2-k}}{F(N)f^{\prime}(N)}\int_{1+\varepsilon}^{\infty} s^{k-1}e^{-\frac{F(N)}{f(N)}s}\,ds
\end{equation*}
\begin{equation*}
+\omega(N)\;\frac{f(N)^{3-k}}{F(N)^{2}f^{\prime}(N)}\int_{1+\varepsilon}^{\infty}s^{k-2}\,e^{-\frac{F(N)}{f(N)}s}\left[1+O\left(\frac{f(N)}{F(N)}\right)\right]ds.
\end{equation*}
Via the scaling $F(N)s=f(N)\,u$ and integration by parts we have
\begin{equation*}
\int_{1+\varepsilon}^{\infty} s^{k-1}\,e^{-\frac{F(N)}{f(N)}s}ds=\left(\frac{f(N)}{F(N)}\right)^{k}\int_{\left(1+\varepsilon\right)\frac{F(N)}{f(N)}}^{\infty}\,\,\,u^{k-1}e^{-u}du
\end{equation*}
\begin{equation*}
=\left(1+\varepsilon\right)^{k-1}\frac{f(N)}{F(N)}e^{-\left(1+\varepsilon\right)\frac{F(N)}{f(N)}}\left[1+O\left(\frac{f(N)}{F(N)}e^{-\frac{F(N)}{f(N)}}\right)\right].
\end{equation*}
Hence, (using again, the definition of $F(\cdot)$ and (\ref{delta}) we get)
\begin{align}
&\sum_{k=0}^{m-1}\frac{F(N)^{k}}{k!}\left\{\int_{1+\varepsilon}^{\infty}s^{k}\left[\int_1^N
f(x)^{-k}\,e^{-\frac{F(N)}{f(x)}s}dx\;\right]ds\right\}\nonumber\\
&=\frac{f(N)^{2}}{F(N)^{2}}\frac{f(N)}{f^{\prime}(N)}e^{-\left(1+\varepsilon\right)\frac{F(N)}{f(N)}}\sum_{k=0}^{m-1}\frac{1}{k!}\left(\frac{F(N)}{f(N)}\right)^{k}\left(1+\varepsilon\right)^{k-1}\left[1+O\left(\frac{f(N)}{F(N)}e^{-\frac{F(N)}{f(N)}}\right)\right]\nonumber\\
&=\frac{\left(1+\varepsilon\right)^{m-2}}{\left(m-1\right)!}\frac{1}{\delta^{m-3}}\,e^{-\varepsilon/\delta}\left(1+O\left(\delta\right)\right).\label{101}
\end{align}
Likewise as $\delta \rightarrow 0^{+}$
\begin{align}
\sum_{k=0}^{m-1}\frac{F(N)^{k}}{k!\,f(N)^{k}}\left[
\int_{1+\varepsilon}^{\infty} s^{k}\,e^{-\frac{F(N)}{f(N)}s}ds\right]=o\left(\frac{1}{\delta^{m-4}}\,e^{-\varepsilon/\delta}\right).\label{102}
\end{align}
In view of (\ref{101}), (\ref{102}), and since $\vartheta \in (0,1)$ is arbitrary, (\ref{b}) implies
\begin{equation*}
I_{3}(N)=\frac{\left(1+\varepsilon\right)^{m-2}}{\left(m-1\right)!}\frac{1}{\delta^{m-3}}\,e^{-\varepsilon/\delta}\left(1+O\left(\delta\right)\right) \quad \text{as}\,\,\, \delta\rightarrow 0^{+}, \,\,\,m=2,3,\dots \,.
\end{equation*}
\hfill $\blacksquare$

\smallskip

\textit{Proof of Lemma 2.3}. We will discuss briefly, the proof for $I_{5}(N)$. The proofs for $I_{4}(N)$ and $I_{6}(N)$ are similar to the
proofs of the results for $I_{1}(N)$ and $I_{3}(N)$ respectively of Lemma 2.2. For $I_{5}(N)$ of (\ref{I5}) and in view of (\ref{AsN}) we have
\begin{equation*}
I_{5}(N) := \int_{1-\varepsilon}^{1+\varepsilon} s\,e^{B_{m}(N;s)}ds.
\end{equation*}
We can treat $I_{5}(N)$ as we treated $I_{2}(N)$ of Lemma 2.2. One gets (as $N \rightarrow \infty$),
\begin{align*}
I_{5}(N)=I_{2}(N)-&
\delta^{2}\int_{\delta^{2-m}\exp\left(-\varepsilon /\delta\right)}^{\delta^{2-m}\exp\left(\varepsilon /\delta\right)}
\ln\left(u\,\delta^{m-2}\right)\nonumber\\
&\times\exp\left\{-\frac{1}{\left(m-1\right)!}\,u\left[1-\delta \ln\left(u\,\delta^{m-2}\right)\right]^{m-2}\right.\nonumber\\
&\left.
-\left(\omega(N)+\frac{1}{\left(m-2\right)!}\right)\,u\,\delta\,\left[1-\delta \ln\left(u\,\delta^{m-2}\right)\right]^{m-3}
\left(1+O\left(\delta\right)\right)\right\}\frac{du}{u},
\end{align*}
We have
\begin{equation}
I_{51}(N)=I_{21}(N)-
\delta^{2}\left(\int_{\delta^{2-m}\exp\left(-\varepsilon /\delta\right)}^{1/\sqrt{\delta}}
+ \int_{1/\sqrt{\delta}}^{\delta^{2-m}\exp\left(\varepsilon /\delta\right)}\right). \label{Aeint5}
\end{equation}
For the second integral of (\ref{Aeint5}) one gets an upper bound (see (\ref{eint6})), namely $O\left(\ln \delta\,\sqrt{\delta}\; e^{-1/\left(m-1\right)!\sqrt{\delta}}\right).$ The first integral of (\ref{Aeint5}) is
\begin{align*}
K_2(\delta):=\int_{\delta^{2-m}\exp\left(-\varepsilon /\delta\right)}^{1/\sqrt{\delta}}
&\ln\left(u\,\delta^{m-2}\right)\exp\left\{-\frac{1}{\left(m-1\right)!}\,u\,\left[1-\delta \ln\left(u\,\delta^{m-2}\right)\right]^{m-2}\right.\nonumber\\
&\left.-\left(\omega(N)+\frac{1}{\left(m-2\right)!}\right)\,u\,\delta\,
\left[1-\delta \ln\left(u\,\delta^{m-2}\right)\right]^{m-3}\right.\nonumber\\
&\left.\times\,\left(1+O\left(\delta\right)\right)\frac{du}{u}\right\}.
\end{align*}
If we treat $K_2(\delta)$ as we treated $K_1(\delta)$ of Lemma 2.2, we can replace the upper limit of the integral $K_{2}(\delta)$ by $\infty$. Thus,
as $\delta\rightarrow 0^{+}$,
\begin{align}
I_{5}(N)=I_{2}(N) -\delta^{2}&\int_{\delta^{2-m}\exp\left(-\varepsilon /\delta\right)}^{\infty}\ln\left(u\,\delta^{m-2}\right)\frac{e^{-u/\left(m-1\right)!}}{u}
\left[1+\frac{m-2}{\left(m-1\right)!}\,u\delta\,\ln\left(u\,\delta^{m-2}\right)\right.\nonumber\\
&\left.-\left(\omega(N)+\frac{1}{\left(m-2\right)!}\right)\,u\,\delta\,\left(1+O\left(\delta\right)\right)
+u^{2}O\left[\delta \ln\left(u\,\delta^{m-2}\right)\right]^{2}\right]du.\label{eint88}
\end{align}
The following asymptotic expansions easy exercises as $x\rightarrow 0^{+}$:
\begin{equation}
\int_{x}^{\infty}\frac{e^{-t}}{t}\ln t\;dt= -\frac{1}{2}\ln^{2} x+\frac{1}{2}\left(\gamma^{2}+\frac{\pi ^{2}}{6}\right)+O\left(x\ln x\right),\label{L}
\end{equation}
\begin{equation}
\int_{x}^{\infty}e^{-t}\ln^{2} t\;{dt}=\left(\gamma^{2}+\frac{\pi ^{2}}{6}\right)-x\ln^{2} x+O\left(x\ln x\right).\label{M}
\end{equation}
Applying (\ref{L}), (\ref{M}) in (\ref{eint88}), and using (\ref{eint11}),  one arrives at
\begin{align*}
I_{5}(N)=&\varepsilon+\frac{\varepsilon^{2}}{2}+\left(m-2\right)\delta\ln \delta+\left[\ln\left(m-1\right)!-\gamma\right]\delta
-\frac{\left(m-2\right)^{2}}{2} \delta^{2}\ln^{2} \delta\nonumber\\
&\,\,\,+\left[\left(m-2\right)^{2}-\left(m-2\right)\left(\ln\left(m-1\right)!-\gamma\right)\right]\delta^{2}\ln \delta\nonumber\\
&\,\,\,+\left[\left(m-2\right)\ln\left(m-1\right)!-\left(m-2\right)\gamma-\omega(N)\left(m-1\right)!-\left(m-1\right)\right.\nonumber\\
&\left.\,\,\,\,\,\,\,\,\,-\frac{1}{2}\left(\ln\left(m-1\right)!\right)^{2}-\frac{1}{2}\left(\gamma^{2}+\frac{\pi^{2}}{6}\right)+\gamma \ln\left(m-1\right)!\right]\delta^{2}+O\left(\delta^{3}\left(\ln\delta\right)^{2}\right).
\end{align*}
\hfill $\blacksquare$

\section{Concluding remarks and comparison with earlier works}
The main task of this paper was to obtain the limiting distribution of the random variable $T_{m}(N)$ (the number of trials a collector needs in order to obtain $m$ complete sets of all $N$ different types of coupons). A key feature in our approach is that for each integer $N > 0$, one can create a probability measure
$\pi _N =\{p_1,...,p_N\}$ on the set of types $\{1,...,N\}$ by taking
\begin{equation*}
p_j = \frac{a_j}{A_N},
\qquad \text{where}\quad
A_N = \sum_{j=1}^N a_j,
\end{equation*}
where $\alpha = \{a_j\}_{j=1}^{\infty}$ is a given sequence of positive numbers.
Under this setup $p_j$ depends on $\alpha $ and $N$. Thus, given $\alpha $, it makes sense to consider the asymptotic behavior of the moments and the variance of the random variable $T_{m}\left(N\right)$. Moreover, since the leading term of $A_N $ is, generally, easy to be found, one focuses in the asymptotics of $E_{m}(N;\alpha)$ and $Q_{m}(N;\alpha)$ (see formulae (\ref{10}), (\ref{13})). Theorem 2.1 separates the problem in classes of growing
(Case I) and decaying sequences (Case II) $\alpha$. For Case I the asymptotics of the expected value and the variance of $T_{m}(N)$ follow easily and \textit{depend} on $m$ (Theorem 2.2).
Having those asymptotics, we were able to establish Theorem 4.1, which gives the limiting distribution of $T_{m}(N)$ (appropriately normalized) as $N \to \infty$ . Notice that the proof of Theorem 4.1 is new, since only for the
case $m = 1$ \cite{DP} the result follows from known theorems together with the asymptotics of the expected value and the variance.\\
In Case II we examine sequences $\alpha = \{a_j\}_{j=1}^{\infty}$ of the form $a_{j}=f(j)^{-1}$, where $f(\cdot)$ satisfies some rather weak conditions (see (\ref{C1a}), (\ref{C1})). In order to apply the general Theorem $N$ (see Subsection 4.2), we need to come up with appropriate
sequences $b_N$ and $k_N$. Here, our asymptotics for $E\left[\, T_{m}\left(N\right)\,\right]$ and $V\left[\, T_{m}\left(N\right)\,\right]$ indicate specifically how to choose $b_N$ and $k_N$ (see Theorem 4.2). Furthermore, in Subsection 4.3 we discuss the case where $f(x)$ does not satisfy the conditions of \eqref{C1}. Formula \eqref{SD1} presented there is completely new. Recall that P. Erd\H{o}s and A. R$\acute{e}$nyi (1961)
\cite{E-R} proved a limit theorem for the case of equal probabilities. Thus, our paper generalizes that result for a large class of coupon probabilities.\\

The computation of the asymptotics of $E\left[\, T_{m}\left(N\right)\,\right]$ and $V\left[\, T_{m}\left(N\right)\,\right]$ in Case II is quite involved. The heart of our analysis rests in Lemma 2.1, which determines the behavior of the quantity $\sum_{j=1}^{N}e^{-\frac{F(N)}{f(j)}}S_{m}(\frac{F(N)}{f(j)}s)$. Hence, it is necessary to rewrite $E_{m}(N;\alpha)$ and $Q_{m}(N;\alpha)$ as in (\ref{b3a}) and (\ref{QM}) respectively. It turns out that we have to compute the fifth asymptotic term of $E_{m}(N;\alpha)$ and the sixth term of $Q_{m}(N;\alpha)$, so that the leading term of the variance emerges (notice that this term is independent of the number $m$). In an earlier work \cite{DP}, the authors established these formulas for $m=1$. The main difference here comes from the limit $\lim_{N}\int_{1}^{N}e^{-\frac{F(N)}{f(x)}s}S_{m}\left(\frac{F(N)}{f(x)}s\right)dx$, (see (\ref{LL2a}) versus (\ref{LL1}) and (\ref{LL2})). This limit causes some difficulties. Its different values (for $m=1$, $m=2$, and $m \geq 3$) explain the reason for considering the formulas (\ref{b3a}) and (\ref{QM}).

\bigskip

\textbf{Acknowledgments.} The authors wish to thank professor Amir Dembo for reading the manuscript and for making valuable comments and
suggestions.\\

\end{document}